\documentclass[11pt]{amsart}

\usepackage{amsmath,amsthm,amsfonts,amssymb,bm,graphicx}
\usepackage[alphabetic,initials]{amsrefs}
\usepackage
{hyperref}
\hypersetup{colorlinks=true,citecolor=blue,linkcolor=blue,urlcolor=blue,
pdfstartview=FitH }

\theoremstyle{plain}
\newtheorem{theorem}{Theorem}[section]
\newtheorem{cor}[theorem]{Corollary}
\newtheorem{lemma}{Lemma}[section]
\newtheorem{prop}{Proposition}[section]

\theoremstyle{definition}

\newtheorem{remark}{Remark}[section]

\renewcommand{\geq}{\geqslant}
\renewcommand{\leq}{\leqslant}

\DeclareMathOperator{\PGL}{PGL}
\DeclareMathOperator{\SL}{SL}

\DeclareMathOperator{\PSO}{PSO}
\def\stacksum#1#2{{\stackrel{{\scriptstyle #1}}
{{\scriptstyle #2}}}}

\newcommand{\sym}{\mathrm{sym}}

\newcommand{\whZ}{\widehat{\Zz}}

\newcommand{\eps}{\varepsilon}

\newcommand{\Cc}{\mathbb{C}}

\newcommand{\Qq}{\mathbb{Q}}

\newcommand{\Zz}{\mathbb{Z}}

\newcommand{\ov}[1]{\overline{#1}}

\newcommand{\mfb}{{\mathfrak{b}}}

 \newcommand{\bfb}{{\mathbf {b}}}

	\newcommand{\mcA}{{\mathcal{A}}}
	\newcommand{\BA}{{\mathbb {A}}}

	 \newcommand{\BH}{{\mathbb {H}}}

	\newcommand{\BQ}{{\mathbb {Q}}} 
	\newcommand{\Rr}{{\mathbb {R}}}
\newcommand{\Zp}{{\mathbb {Z}_p}}

\newcommand{\peter}[1]{\langle{#1}\rangle}

	\newcommand{\Hh}{\BH}
		
	\newcommand{\Qp}{\BQ_p}
	\newcommand{\Qv}{\BQ_v}
	\newcommand{\Aa}{\BA}

	\newcommand{\bash}{\backslash}

\newcommand{\ra}{\rightarrow}
\newcommand{\mods}[1]{\,(\mathrm{mod}\,{#1})}
\def\peter#1{\langle #1\rangle}

\newcommand{\bfn}{\mathbf{n}}
\newcommand{\rmG}{\mathrm{G}}

\newcommand{\GQS}{\mathrm{G}(\Qq_S)}

\newcommand{\vphi}{\varphi}

\numberwithin{equation}{section}

\makeatletter
\def\Ddots{\mathinner{\mkern1mu\raise\p@
\vbox{\kern7\p@\hbox{.}}\mkern2mu
\raise4\p@\hbox{.}\mkern2mu\raise7\p@\hbox{.}\mkern1mu}}
\makeatother
\makeatletter
\DeclareRobustCommand\widecheck[1]{{\mathpalette\@widecheck{#1}}}
\def\@widecheck#1#2{%
    \setbox\z@\hbox{\m@th$#1#2$}%
    \setbox\tw@\hbox{\m@th$#1%
       \widehat{%
          \vrule\@width\z@\@height\ht\z@
          \vrule\@height\z@\@width\wd\z@}$}%
    \dp\tw@-\ht\z@
    \@tempdima\ht\z@ \advance\@tempdima2\ht\tw@ \divide\@tempdima\thr@@
    \setbox\tw@\hbox{%
       \raise\@tempdima\hbox{\scalebox{1}[-1]{\lower\@tempdima\box
\tw@}}}%
    {\ooalign{\box\tw@ \cr \box\z@}}}
\makeatother

\begin{document}

\author{Valentin Blomer}
\address{Mathematisches Institut, Endenicher Allee 60, 53115 Bonn, Germany} \email{blomer@math.uni-bonn.de}
 
\author{Philippe Michel}
\address{EPFL-SB-MATH-TAN, Station 8, 
1015 Lausanne, Switzerland}\email{philippe.michel@epfl.ch}

 \title{The unipotent mixing conjecture}

\thanks{First author supported in part by Germany’s
Excellence Strategy grant EXC-2047/1 - 390685813 and ERC Advanced Grant 101054336. Second author partially supported by the SNF grant 200021\_197045.}


\begin{abstract}  We show that shifted pairs of discrete or continuous low-lying horocycles equidistribute in the product space   of two modular curves. 
 \end{abstract}


\setcounter{tocdepth}{2}  \maketitle 

\maketitle

\begin{center} \emph{Dedicated to Peter Sarnak with admiration}\end{center}


\section{Introduction}

\subsection{Main results}

A classical result of Sarnak \cite{Sa} says in a special case that the image of a low-lying horocycle $$\{x + i/T \mid x \in [0, 1]\}$$ of height $1/T$ equidistributes on the modular curve $X = {\rm SL}_2(\mathbb{Z}) \backslash \mathbb{H}$ with respect to the usual hyperbolic probability measure $$d\mu_X(z)=\frac{3}{\pi} \frac{dx\, dy}{y^2}$$ as $T \rightarrow \infty$. 

Similarly, one can consider  for a prime $q$, say,  ``discrete'' low-lying horocycles 
$$H_q:=\Big\{\frac{a+i}{q} \mid a\, (\text{mod } q)\Big\}.$$
  Except for the point $iq \in \mathbb{H}$, this is the image of $\SL_2(\Zz)i$ by the Hecke correspondence $T_q$ (which has degree $q+1$), and  non-trivial bounds for eigenvalues of Hecke operators imply that these discrete low-lying horocycles also equidistribute on $X$ as $q \rightarrow \infty$. 
  
  It is natural to investigate whether equidistribution persists when some additional constraints are imposed on the entries $a$. Identifying the set of $a$'s with integers contained in the interval $[0,q-1]$, one can restrict $a$ to belong to a subinterval $I\subset [0,q-1]$ (see \cite{Strom}) or vary $a$ along a subsequence of integers of a special shape like $a=[n^c],\ c=1.1$ \cite{Ven} or ask for $a$ to vary along the primes.
  
  This last case was studied by Sarnak and Ubis \cite{SU} who proved, assuming the Ramanujan-Petersson-Selberg conjecture that, as $q\ra\infty$, the ``prime'' $q$-Hecke points 
  $$H_q^{\mathrm p}:=\Big\{\frac{p+i}{q},\ 1\leq p\leq q-1,\ p\hbox{ prime}\Big\}.$$
  become dense in $X$; more precisely  any weak limit, $\mu^{\mathrm{p}}$ say, of the uniform probability measures supported by the $H_q^{\mathrm p}$ satisfies
  $$\frac15\mu_X\leq \mu^{\mathrm p}\leq \frac95\mu_X.$$
  As for the proof, the combinatorial decomposition of the characteristic function of the primes lead naturally to handling sums of ``Type I'' and ``Type II''. The treatment of type I sums is closely related to restricting  $a$ to a sub-interval while the treatment of type II leads a new equidistribution problem: the joint equidistribution  of sets of the shape
  $$H_{q,(b_1,b_2)}(N):=\Big\{\Big(\frac{b_1a+i}{q},\frac{b_2a+i}{q}\Big) \mid  1\leq a\leq N\Big\}\subset X\times X$$
  on the product $X\times X$; here $0< b_1<b_2$ are distinct integers bounded by a fixed power of $q$ and with $N$ satisfying $b_1N,b_2N\leq q$.
  
 In the present paper we investigate the distribution properties of the above sets when we let $a$ vary over the {\em complete set} of congruence classes modulo $q$ and look at how
$$H_{q,b}:=\Big\{\Big(\frac{a + i}{q}, \frac{ab  + i}{q}\Big) \mid a\mods q\Big\}$$
distributes in  the product $X \times X$ as $q \rightarrow \infty$ as we allow $b$ to possibly vary with $q$. 
Clearly  if $b=1$ we cannot expect full equidistribution to hold: $H_{q, b}$ is simply trapped in the diagonally embedded copy of $X$ in $X\times X$ (and equidistributes along it). A  similar phenomenon will occur for any fixed (i.e.\ independent of $q$) residue class $b$. One might optimistically conjecture that something like $q \| b/q \| \rightarrow \infty$ (where $\| . \|$ denotes the distance in $\mathbb{R}/\mathbb{Z}$) may suffice to ensure equidistribution, but the situation is more subtle as we shall see in a moment. 

 It turns out to be useful to  introduce the lattice
\begin{equation}\label{lattice}
\Lambda_{q,b}=\{(n_1,n_2)\in \Zz^2,\ n_1+n_2b\equiv 0\mods q\}\subseteq \Zz^2.
\end{equation}
This lattice has covolume $q$ and we denote by $$s(q;b)=\min_{\mathbf{0}\not=\bfn\in \Lambda_{q,b}}\|\bfn\|$$
its minimum (i.e.\ the minimal euclidean norm of a non-zero element). 
It is well-known from the geometry of numbers that 
\begin{equation}\label{srange}
	s(q;b) \ll q^{1/2}. 
\end{equation}
In particular $s(q;b)\ra\infty$ implies $ q\ra\infty.$ We are now ready to state our first main result.

\begin{theorem}\label{main1}  Let $q$ be a large prime and $b \in (\mathbb{Z}/q\mathbb{Z})^{\times}$. 
Then for pairs $(q,b)$ such that $s(q;b)\rightarrow \infty$ the set
\begin{equation}
	\label{horocyclepairs}
	H_{q,b}=\Big\{\Big(\frac{a + i}{q}, \frac{ab  + i}{q}\Big) \mid a\mods q\Big\} \subseteq X \times X
\end{equation}
becomes equidistributed with respect to the product of hyperbolic measures $\mu_{X\times X}=\mu_X\otimes\mu_X.$ 
\end{theorem}
Once the case of two factors is available, it is possible to obtain equidistribution for more (\cite{EL}*{Corollary 1.5}):  
\begin{cor}
    \label{maincor}
    For $d\geq 2$ and a $d$-tuple of congruence classes $\bfb=(b_1,\cdots,b_d)\in ((\Zz/q\Zz)^\times)^d$
    let $$s(q;\bfb):=\min_{i\not=j}s(q;b_ib_j^{-1}).$$
    Then for pairs $(q,\bfb)$ such that $s(q,\bfb)\rightarrow \infty$ the set
    $$H_{q,\bfb}:=\Big\{\Big(\frac{ab_1 + i}{q},\cdots, \frac{ab_d  + i}{q}\Big) \mid a\mods q\Big\}\subseteq X^d$$
    becomes equidistributed with respect to the product of hyperbolic measures $\mu_{X^d}=\mu^{\otimes d}_X.$
\end{cor}
\begin{remark} Theorem \ref{main1} and Corollary \ref{maincor} hold true with $X$ replaced by any fixed modular curve  $X_0(D)=\Gamma_0(D)\bash \Hh$ for $D\geq 1$ and $$\Gamma_0(D)=\Big\{\begin{pmatrix}
    a&b\\c&d
\end{pmatrix}\in\SL_2(\Zz), \ c\equiv 0\mods D \Big\}.$$
 In fact, for the proof we will need such an extension to a modular curve of higher level.
\end{remark}
\begin{remark}
If we denote 
for a continuous, compactly supported function $\vphi:X\times X\ra \Cc$ by
\begin{equation}\label{Weylsum0}
\mathcal{W}_\vphi(b;q):=\frac{1}q\sum_{a\mods q}\vphi\Big(\frac{a + i}{q}, \frac{ab  + i}{q}\Big)
\end{equation}
 the corresponding Weyl sum, then one approach to Theorem \ref{main1} is to show that 
\begin{equation}
	\label{Weylsum}
	\mathcal{W}_\vphi(b;q)=\int_{X\times X}\vphi(z_1,z_2)d\mu_{X\times X}(z_1,z_2)+o_\vphi(1). 
\end{equation}
If we assume the Ramanujan-Petersson conjecture for the Hecke eigenvalues of Maa{\ss} cusp forms for the group ${\rm SL}_2(\mathbb{Z})$ and $\vphi$ is cuspidal (in both components), we prove \eqref{Weylsum} with a polylogarithmic rate in the form 
$$o_\vphi(1)\ll_{\vphi,\eps} s(q;b)^{-1+\varepsilon} + \log^{-\delta}p$$
for some fixed $\delta>0$ and any $\varepsilon > 0$. See Section \ref{sec12} for more details. 
\end{remark}

\begin{remark}
	The condition   $s(q;b) \rightarrow \infty$ means that $b/q$ must not be too close to a rational number with fixed denominator (like $2/5$), so $b/q$ must be, in a rather weak sense, ``badly approximable''. We indicate in Remark \ref{sisoptimal} below why this condition  is probably sharp.  
 \end{remark}

 \begin{remark}
 Theorem \ref{main1} can probably be generalized to composite numbers $q$ at the cost of slightly more work. For instance, some care must be taken if $q$ has many small prime factors which might collide with the auxiliary primes to be chosen later. 
\end{remark}

We prove a similar theorem for ``continuous'' low-lying pairs of horocycles.

 \begin{theorem}\label{main2}  
 Let $T > 1$, $y \in [1, 2]$ and write $y = a/q + O(1/qQ)$ for positive coprime integers $a, q$ with $q \leq Q := T^{0.99}$. Let $I \subseteq \mathbb{R}$ be a fixed non-empty interval. Then 
$$\Big\{\Big(\frac{x + i}{T}, \frac{xy + i}{T}\Big) \mid x \in I\Big\} \subseteq X \times X$$
equidistributes as $T \rightarrow \infty$ for pairs $(y, T)$ with $q \rightarrow \infty$. 
 \end{theorem}
 
 The condition on $y$ is similar: $y$ must be in some weak sense badly approximable by rational numbers with small denominator. There is some flexibility in the definition of $Q$ as the proof shows. Also Theorem \ref{main2} can be generalized to $d$ factors in a similar way. \\ 
 
 These results are new examples of a growing family of equidistribution results inspired by Sarnak's classical theorem \cite{Sa} mentioned above, as well as by Duke's celebrated equidistribution theorems for Heegner points and closed geodesics with increasing discriminants on modular curves \cite{Du}; the latter needs crucially subconvexity for twisted $L$-functions \cite{DFI} as an input.

A new generation of equidistribution problems was posed by the second author and Venkatesh in their ICM address \cite{MV} in the context of Duke's theorem. Let  $H(D)$ be the set of Heegner points in the imaginary quadratic field $K = \mathbb{Q}(\sqrt{-D})$. Let $[\mfb]$ be an element in the class group of $K$ (which acts on $H(D)$), and consider the set of pairs $$\{(z, [\mfb]\star z) \mid z\in H(D)\}\in X\times X.$$ As $D \rightarrow \infty$, does this equidistribute  in the product space $X \times X$?  Similarly to the assumption   $s(q;b) \rightarrow \infty$ in Theorem \ref{main1}, a necessary condition  to escape from the diagonal is that the norm of the smallest ideal in $[\mfb]$ tends to infinity, and it is reasonable to conjecture that this is sufficient. This is (a special case of) the \emph{mixing conjecture}. This conjecture is still open but two different conditional proofs were given by Khayutin  \cite{Kh} and Blomer-Brumley-Khayutin \cite{BBK}. One of the key ingredients in the former  is a deep result in homogeneous dynamics due to Einsiedler and Lindenstrauss \cite{EL} classifying under suitable condition joinings in products of locally homogeneous spaces. As we will see below this result also play a crucial role in the present paper. Further variations on this theme can be found in \cites{AES,AEW,ALMW,BB}, and we recommend in particular the Bourbaki seminar by M.\ Aka \cite{AkaBBKI}.


\begin{remark} The astute reader might have  noticed that $H_q$  is essentially the set $H(-4q^2)$ of Heegner points of the highly non-maximal order of conductor $-4q^2$. 

More precisely, if $q \equiv 3$ (mod 4), we have $ \# H(-4q^2) = \frac{1}{2}(q+1)$, and $H_q$ corresponds to the $q$   primitive quadratic forms $F_{q, a} := (q^2, 2qa, a^2 + 1)$ for $a\,  (\text{mod } q)$.  
The forms $F_{q, a_1}$ and $F_{q, a_2}$ are equivalent if and only if $a_1 a_2 \equiv - 1$ (mod $q$) or $a_1 \equiv a_2$ (mod $q$), so that   $H_q$ covers $H(-4q^2)$   twice except for the form $F_{q, 0}$ which corresponds uniquely to the identity in $H(-4q^2)$. 

If $q \equiv 1$ (mod 4), then  $\# H(-4q^2) = \frac{1}{2}(q-1)$, and the same analysis holds, except that the two forms $(q^2, 2qa, a^2 + 1)$ with $a^2 \equiv -1$ (mod $q$) are not primitive any more, and the corresponding two representatives of the same point in ${\rm SL}_2(\mathbb{Z}) \backslash \mathbb{H}$ are not counted by $H(-4q^2)$. 

However the group action in $H(-4q^2)$ is quite different from the multiplication in $H_q$. \\

We recall that Khayutin's work \cite{Kh}  is tailored for fundamental discriminants and uses quite crucially the fact that the underlying order is maximal (or at least not far away from being maximal, cf.\ \cite{Kh}*{Section 1.8.2}). 
While our proof also builds crucially on the work of Einsiedler-Lindenstrauss, our approach is different from Khayutin's: instead of bounding a correlation function between two measures, we prove, by direct arguments, that the Weyl sums attached to a well chosen cuspidal test function $\vphi$ converge to the correct limit (i.e.\ $0$); the classification theorem of Einsiedler-Lindenstrauss then allows us to bootstrap this partial information to full equidistribution.
\end{remark}

\begin{remark} Unlike other examples, our result  does not require any splitting conditions: in our case such conditions are in fact automatically satisfied. As such, this seems to be one of the first \emph{unconditional} instances of a mixing type conjecture (see also \cite{Saw} as well as \cite{Tsi} for proofs of function field versions of this conjecture). 
\end{remark}

\begin{remark} 
Lindenstrauss, Mohammadi and Wang \cite{LMW} have obtained very general and effective forms of equidistribution  with polynomial decay rates for certain types of one parameter unipotent flows on products of modular curves. It is not clear to us whether this would cover all cases of Theorem \ref{main2} (as the analysis in \cite{LMW} depends on the injectivity radius of the base point $x_0$) and then whether this could be transferred to the discrete  case treated in Theorem \ref{main1}. 
\end{remark}

\subsection{Principle of the proof}\label{sec12}

By the spectral decomposition and Weyl's equidistribution criterion, it is sufficient to bound non-trivially the Weyl sum $\mathcal{W}_\vphi(b;q)$  for $\vphi=f_1\otimes f_2$ when $f_1$ and $f_2$ are either constant functions or (non-constant) Hecke eigenforms. The case where $f_1$ or $f_2$ is constant follows from the  equidistribution of \eqref{horocyclepairs} when projected to each factor $X$. To handle the remaining sums and prove  Theorem \ref{main1} we 
use different techniques depending on how fast $s(q;b)$ approaches $\infty$ within the range \eqref{srange}, a principle that is familiar from \cite{EMV}, \cite{Kh} or \cite{BBK}.

More precisely, let $\theta=7/64$ be the best known approximation towards the Ramanujan-Petersson conjecture \cite{KimSar} and suppose we are in the range
$$s(q;b)\leq q^{1/2-2\theta-\eta}$$
for some fixed $\eta>0$.  Bounds for the Weyl sums can be achieved by resolving some versions of the shifted convolution problem (see Proposition \ref{thm1}). In this regime, the Ramanujan-Petersson conjecture is not necessary and one can show that equidistribution holds with a polynomial decay rate $o_\vphi(1)\ll q^{-\delta(\eta)}$ in the notation of \eqref{Weylsum}. 

In the remaining range $$q^{1/2-2\theta-\eta}\leq s(q;b)\ll q^{1/2},$$ methods from harmonic analysis do not work so well at the moment. Using sieve methods instead together with the analytic theory of multiplicative functions, we are able to bound non-trivially the Weyl sums 
\eqref{Weylsum} when $f_1$ and $f_2$ are   cuspidal and both satisfy the Ramanujan-Petersson conjecture. The Ramanujan-Petersson conjecture is needed to implement the sieving argument as it insures  that the (multiplicative) Hecke eigenvalue functions $n\mapsto \lambda_{f_i}(n)$ are bounded in absolute value by the divisor function. This reduces the problem to that of obtaining a sufficiently good bound for  sums of the Hecke eigenvalues along the  primes, namely 
\begin{equation}\label{goodSTbound}
	\sum_\stacksum{p\leq z}{p\, \text{prime}}\frac{|\lambda_{f_1}(p)|+|\lambda_{f_2}(p)|}p\leq (2-\delta) \log\log z 
\end{equation}
for some fixed $\delta>0$ and sufficiently large $z$. Such a bound, which follows from suitable approximations to the Sato-Tate conjecture for cusp forms, yields a power saving in $\log p$ on the size of the Weyl sum (in a way similar to \cite{Ho}). Unfortunately, if $f_1$ or $f_2$ is an Eisenstein series, the distribution properties for their Hecke eigenvalues  do not allow to obtain \eqref{goodSTbound}. Nevertheless it will be useful to remember that if $f_1$ and $f_2$ are both {\em CM forms} (cusp forms attached to Hecke characters of   quadratic fields) the Ramanujan-Petersson conjecture holds and \eqref{goodSTbound} is unconditional.

To deal with the rest of the spectrum  we take a very different route and make use of  a powerful measure classification theorem of Einsiedler-Lindenstrauss \cite{EL}*{Thm 1.4}. Let $\mu_{q,b}$ denote the uniform probability measure on $X\times X$ supported by $H_{q,b}$. As we explain below in Section \ref{sec5}, $\mu_{q,b}$ is the projection to $X\times X$ of a measure $\mu^S_{q,b}$ on a suitable $S$-arithmetic quotient $$\Gamma_S\times\Gamma_S\bash \rmG(\Qq_S)\times\rmG(\Qq_S),\ \rmG=\PGL_2$$ and (by non-trivial bounds for Hecke eigenvalues in the discrete case or Sarnak's theorem in the continuous case) any weak-$\star$ limit of the $\mu^S_{q,b}$ projects to the Haar measure on each factor, i.e.\ it defines a joining of the Haar measures. Moreover, these limits are invariant under the action of a rank $2$ diagonalisable subgroup, namely the diagonal subgroup of $\rmG(\Qq_S)\times\rmG(\Qq_S)$ generated by 
$$(t_1,t_1)\hbox{ and }(t_2,t_2)\hbox{ where }t_i=\begin{pmatrix}
    q^{-1}_i&\\&q_i
\end{pmatrix},\ i=1,2$$
for  $q_1,q_2$ two primes distinct from $q$. This is an immediate consequence of the invariance of $H_{q,b}$ under multiplicative shifts that we have already observed:
$$
	H_{q,b} =\Big\{\Big(\frac{aq_i^2 + i}{q}, \frac{aq_i^2b  + i}{q}\Big) \mid a\mods q\Big\},\quad  i=1,2.$$
	
	By the measure classification theorem, any such joining is a convex combination of {\em algebraic} measures which in the present situation are either the (image of the) Haar measure on the full product space or the Haar measure along diagonal $\rmG$-orbits in the product. We will have to exclude the latter possibility. This will follow -- in a rather subtle way -- from what we have already proved. Testing the measure against a carefully selected test function $\vphi = f_1 \otimes f_2$, where both $f_1$ and $f_2$ are suitable vectors in a certain CM automorphic representation (which satisfies the Ramanujan conjecture), we show that the corresponding Weyl sum  would have a \emph{positive} limit if the measure were not the full product measure.

	\begin{remark}
		This principle of  using measure rigidity  to bootstrap the  evaluation of Weyl sums for a small portion of the spectrum  to full a equidistribution statement has occured in the past notably in \cite{ELMV}  (for the Siegel Eisenstein series for $\SL_3$). This principle is likely to be used again in higher rank situations in connection to functoriality: for instance we hope that the results of this paper will be useful for equidistribution problems associated with ${\rm GSp}_4$ (on using Saito-Kurokawa lifts). 	
	\end{remark}

\begin{remark} Notice that the two entries $(a\mods q,ab\mods q)$ in Theorem \ref{main1} are related by a linear equation. Another interesting question occurring in classical problems from analytic number theory is whether joint equidistribution holds for tuples of discrete horocycles whose entries are related by more general algebraic equations. For instance,  the second named author together with Einsiedler and Lindenstrauss (\cite{ELM}) established joint equidistribution in the case of monomial relations: given two fixed integers $1\leq k<l$, then as $q\ra\infty$ amongst the primes and for any $b\in(\Zz/q\Zz)^\times$, the set
$$H_{k,l,q,b}:=\Big\{\Big(\frac{a^k + i}{q}, \frac{ba^l  + i}{q}\Big) \mid a\mods q\Big\}$$
 becomes equidistributed on $X\times X$. The proof uses again crucially the measure classification results of \cite{EL}, but the application is {\em simpler}: the sets $H_{k,l,q,b}$ are then invariant under the rank $2$ subgroup generated by $(t_1^k,t_1^l)$ and $(t_2^k,t_2^l)$, and the fact that $1\leq k<l$ excludes the possibility of having measures supported along diagonal $\rmG$-orbits. The combination of the present results with those of \cite{ELM} will be the topic of a future work, joint with Einsiedler and Lindenstrauss.
\end{remark}

\section{Shifted convolution problems}

In this section we prepare the stage for a proof of Theorem \ref{main1} in the ``low regime'' case, i.e.\ when
$$s(q;b)\leq q^{1/2-2\theta-\delta}$$
for some $\delta>0$. 

We start with a variation of  \cite{Bl}*{Theorem 1.3} with an additional summation over $h$.

\begin{prop}\label{thm1} Let $N_1N_2, d, l_1, l_2 \in \mathbb{N}$ with $(N_1N_2l_1l_2, d) = 1$, $H, M_1, M_2, P_1, P_2 \geq 1$. Let $f, g$ be two (holomorphic or Maa{\ss}) cuspidal newforms of levels $N_1, N_2$ respectively and central characters $\chi_f$, $\chi_g$ with Hecke eigenvalues $\lambda_f(m), \lambda_g(m)$. Let $G$ be a smooth function supported on $[M_1, 2M_1] \times [M_2, 2M_2]$ satisfying $\| G^{(i, j)}  \|_{\infty} \ll_{i, j} (P_1/M_1)^i (P_2/M_2)^j$ for $i, j \in \mathbb{N}_0$. For $H \leq h \leq 2H$ let $|\alpha(h)| \leq 1$. Then\\
\begin{displaymath}
\begin{split}
\mathcal{D} &:= \sum_{\substack{H \leq h \leq 2H\\ d \mid h, (d, h/d) = 1}} \alpha(h) \Big|\sum_{l_1m_1 \pm l_2m_2 = h} \lambda_f(m_1)\lambda_g(m_2) G(m_1, m_2)\Big| \\
&\ll  (l_1M_1 + l_2M_2)^{1/2 + \theta + \varepsilon} \Big(\frac{H}{d}\Big)^{1/2} \Big(1 + \frac{H/d}{l_1l_2  (1 +  H/l_2M_2)}\Big)^{1/2}
\end{split}
\end{displaymath}
with an implied constant depending on $\varepsilon$ and polynomially on $P_1, P_2$ and the conductors of $f$ and $g$. \end{prop}

The proof depends on the following  bound for averages of twisted Kloosterman sums.   Here we work unconditionally and denote by $\theta \leq 7/64$ an admissible exponent towards the Ramanujan-Petersson conjecture.  

\begin{lemma} \label{prop} Let $P_0, P_1, P_2, S, H, Q \geq 1$, $d, N \in \mathbb{N}$ with $(d, N) = 1$. Let $\chi$ be a (possibly trivial) Dirichlet character modulo $N$. Let $u$ be a smooth function with support in $[H, 2H] \times [S, 2S] \times  [Q, 2Q]$ satisfying  $\| u^{(ijk)}\|_{\infty} \ll (P_0/H)^{i} (P_1/S)^j(P_2/Q)^k$ for $0 \leq i, j, k \leq 2$. Let $a(h)$ and $b(s)$ for $H \leq h \leq 2H$, $S \leq s \leq 2S$ be sequences of complex numbers. 
Then
\begin{displaymath}
\begin{split}
\sum_s\sum_{N \mid q} \sum_{\substack{d \mid h\\ (d, h/d) = 1}} &a(h) b(s) S_{\chi}(\pm h, s, q)  u(h, s, q) \ll  Q \Big( \sum_{s} |b(s)|^2\Big)^{1/2} \Big(1 + \frac{HS}{Q^2} + \frac{S}{N}\Big)^{1/2} \\
& \times d^{\theta}\Big( \sum_{d\mid h} |a(h)|^2\Big)^{1/2} \Big(1 + \frac{H/d}{N  (1 +  HS/Q^2)}\Big)^{1/2} \Big(1 + \Big(\frac{HS}{Q}\Big)^{-\theta}\Big)(NHSQ)^{\varepsilon}
\end{split}
\end{displaymath}
with an implied constant that depends on $\varepsilon$ and polynomially on $P_0, P_1, P_2$.\end{lemma}

\textbf{Proof.} We follow the proof of \cite{Bl}*{Proposition 3.5}. Instead of $U_h(t, q)$ we work with
$$U(t_0, t_1; q) = \int_{-\infty}^{\infty}  \int_{-\infty}^{\infty}  u\Big(x, y, \frac{4\pi \sqrt{xy}}{q}\Big) \frac{4\pi \sqrt{xy}}{q} e(-t_0x - t_1y) \, dx\, dy$$
which satisfies
$$\frac{\partial^n}{\partial q^n} U(t_0, t_1; q) \ll_n \Big(1 + \frac{H|t_0|}{P_0}\Big)^{-2}\Big(1 + \frac{S|t_1|}{P_1}\Big)^{-2}HSQ \Big( \frac{QP_2}{\sqrt{HS}}\Big)^n$$
and
$$q u(h, s, q) =  \int_{-\infty}^{\infty}  \int_{-\infty}^{\infty}  U\Big(t_0, t_1;\frac{4\pi \sqrt{hs}}{q}\Big) e(t_0h + t_1s) \, dt_0\, dt_1.$$
As in \cite{Bl}*{Proposition 3.5} we now apply the Kuznetsov formula, Cauchy-Schwarz and the spectral large sieve for cusp forms of level $N$ and character $\chi$. Not that, as $(d, N) = (d, h/d) = 1$, we can extract the divisibility condition $d\mid h$ at the cost of a factor $d^{\theta}$ in the $h$-sum and reduce the length of the $h$-sum to $H/d$. 
 We then arrive in the same way at \cite{Bl}*{(3.19)} except that the factor $Th^{\theta}$ is replaced with $$d^{\theta}\Big(T^2 + \frac{H/d}{N}\Big)^{1/2} \Big(\sum_{d \mid h} |a(h)|^2\Big)^{1/2},$$
and the result follows with the same choice of $T$ as in \cite{Bl}*{(3.20)}. \qed\\

For the \textbf{proof of Proposition \ref{thm1}} we follow literally the proof of \cite{Bl}*{Theorem 1.3} in \cite{Bl}*{Section 4} with the only minor modification that $f$ and $g$ may be two different cusp forms of potentially levels, they may have non-trivial central characters, and that the summation condition can have either sign. 
We keep the extra sum over $h$ outside until the very end when we apply our Lemma \ref{prop} with $S = Q^2/l_2M_2$ and  $N = N_1N_2l_1l_2$ as a replacement for \cite{Bl}*{Proposition 3.5}. This replaces the factor $h^{\theta}$ in \cite{Bl}*{(4.19)} with 
$$d^{\theta} \Big(\frac{H}{d}\Big)^{1/2}  \Big(1 + \frac{H/d}{l_1l_2  (1 +  H/l_2M_2)}\Big)^{1/2}, $$
so that 
$$\mathcal{D}  \ll (l_1M_1 + l_2M_2)^{1/2 + \theta + \varepsilon} \Big(\frac{H}{d}\Big)^{1/2} \Big(1 + \frac{H/d}{l_1l_2  (1 +  H/l_2M_2)}\Big)^{1/2}$$
as desired.  \qed
 
\section{Application of a sieve}

  For $X, Y \gg 1$ and $q \in \mathbb{N}$ we denote by $\mathcal{C}_q(X, Y)$ the set of  subsets   $\mathcal{S} \subseteq \mathbb{N}^2$ contained in a ball of radius $X^{100}$ about the origin and satisfying  
$$\sum_{\substack{\bfn=(n_1,n_2) \in \mathcal{S}\\ d_1 \mid n_1, d_2 \mid n_2}} 1 = \frac{X}{d_1d_2} + O(Y)$$
for all $(d_1d_2, q) = 1$. A typical example is  $\mathcal{S} = \Lambda \cap \mathcal{R}$ where $\Lambda \subseteq \mathbb{Z}^2$ is a sublattice of  index $q$ and shortest nonzero vector of length $s$, $\mathcal{R}\subseteq \mathbb{R}_{>0}^2$ a  non-empty, bounded, connected and simply-connected set with piecewise smooth boundary $\partial \mathcal{R}$, with $$X = \frac{ \text{vol}(\mathcal{R})}{q}, \quad Y = \frac{\text{length}(\partial \mathcal{R})}{s}.$$ Note that for $(d_1d_2, q) = 1$, the pairs $(n_1, n_2) \in \Lambda$ with $d_j \mid n_j$ for $j = 1, 2$ form a sublattice of index $d_1d_2$.

\begin{prop}\label{thm-sieve} Let $\mathcal{S} \in \mathcal{C}_q(X, Y)$. 
Let $\lambda_1, \lambda_2$ be two non-negative multiplicative functions satisfying $\lambda_j(n)  \leq \tau_k(n)$ for some   $k \in \mathbb{N}$. Fix $0 < \gamma < 1/2$ and let $ z = X^{\gamma}$. Then 
$$\sum_{\substack{ (n_1, n_2) \in \mathcal{S}\\ (n_1n_2, q) = 1}} \lambda_1(n_1) \lambda_2(n_2)  \ll_{k,  \gamma,\varepsilon} \frac{X}{(\log z)^2}
\exp\Big( \sum_{p \leq z} \frac{\lambda_1(p) + \lambda_2(p)}{p} \Big) 
+  \frac{X^{1+\varepsilon}}{ z^{1/4}} + X^{\varepsilon} Y z^3$$
for any $\varepsilon > 0$. 
\end{prop}

For $z\geq1 $ and $q \in \mathbb{N}$ let $\mathcal{P}_{q, z}$ be the set of primes $p \leq z$, $p \nmid q$. For $\mathcal{S} \in \mathcal{C}_q(X, Y)$,  $a_1, a_2 \in \mathbb{N}$ and  $y_1, y_2 \geq 1$ 
define 
$$S(\textbf{a}, \textbf{y}) := \sum_{\substack{\bfn=(n_1,n_2)\in \mathcal{S} \\ a_1 \mid n_1, a_2 \mid n_2\\ (n_1/a_1, \mathcal{P}_{q, y_1}) = (n_2/a_2, \mathcal{P}_{q, y_2}) = 1}} 1.$$
As usual we write $P^+(n)$ for the largest prime factor of $n$ (with the convention $P^+(1) = 1$). From a standard sieve, e.g.\ Selberg's sieve, we obtain

\begin{lemma}\label{lem3} For $\mathcal{S} \in \mathcal{C}_q(X, Y)$, 
$ y_1, y_2 \geq 1$, $a_1, a_2 \in \mathbb{N}$ with $(q, a_1 a_2) = 1$ we have
$$S(\textbf{a},  \textbf{y})  \ll   \frac{X}{ a_1a_2 \log (1+y_1) \log (1+y_2)}   +    Y y_1^2y_2^2 . $$
\end{lemma}

\bigskip

For the \textbf{proof of Proposition \ref{thm-sieve}} we decompose a general $n_1$ appearing in the sum as   $$n_1 = p_1^{e_1} \cdots p_{k_1}^{e_{k_1}}  \cdot p_{k_1+1}^{e_{k_1+1}} \cdots p_r^{e_r} = a_1b_1$$ 
with $$p_1 < p_2 < \ldots,  \quad p_j \nmid q$$
and similarly for $n_2 = a_2b_2$ where $k_j$ is maximal with $a_j \leq z = X^{\gamma}$. Let $$\xi :=  (\log z)( \log\log z).$$ We distinguish the following four cases  for $n_1$:
\begin{displaymath}
\begin{split}
&{\rm (I)}\quad   p_{k_1 + 1} \geq z^{1/2}, \quad {\rm (II)} \quad p_{k_1 + 1} < z^{1/2}, a_1 \leq z^{1/2},\\
& {\rm (III)}\quad  p_{k_1 + 1} < \xi, a_1 > z^{1/2}, \quad {\rm (IV)}\quad  \xi \leq p_{k_1 + 1} \leq z^{1/2}, a_1 > z^{1/2}
\end{split}
\end{displaymath}
and similarly for $n_2$. This set-up goes back originally to Erd\"os \cite{Er} and was refined by Wolke \cite{Wo}, Nair-Tenenbaum \cite{NT}, Khayutin \cite{Kh} and others.

1) Consider first pairs  $ (n_1, n_2)$ such that $n_1$ is in case (II).  Then $p_{k_1+1}^{e_{k_1+1}} \geq z^{1/2}$, so $n_1$ is divisible by a prime power $p^e \geq z^{1/2}$ with $p < z^{1/2}$. Let $e_0 = e_0(p, z) \geq 2$ be the smallest positive integer with $p^{e_0} \geq z^{1/2}$.  Estimating $\lambda_1(n_1)\lambda_2( n_2)$ trivially, by Lemma \ref{lem3} the contribution of such pairs is at most
\begin{equation}\label{case1}
\begin{split}
&\ll X^{\varepsilon} \sum_{\substack{p \leq z^{1/2}\\ p \nmid q}} S(p^{e_0}, 1, 1, 1)  \ll X^{1+\varepsilon}\Big(\sum_{p \leq z^{1/4}} \frac{1}{z^{1/2}} + \sum_{z^{1/4} \leq p \leq z^{1/2}} \frac{1}{p^2}\Big) +  z^{1/2} X^{\varepsilon}Y\\
& \ll \frac{X^{1+\varepsilon}}{  z^{1/4}}  +  z^{1/2}X^{\varepsilon}Y. 
\end{split}
\end{equation}
The same argument works with exchanged indices if $n_2$ is in case (II).

2) Next consider pairs  $ (n_1, n_2)$ such that $n_1$ is in case (III). Then $n_1$ is divisible by a number $z^{1/2} < a_1 \leq z$ such that $P^+(a_1) < \xi$. Since
$$\sum_{\substack{ a \leq x\\  P^+(a) \leq \log x \log \log x}}1 \ll x^{\varepsilon},$$
estimating $\lambda_1(n_1)\lambda_2( n_2)$ trivially, the contribution of such pairs is by Lemma \ref{lem3} and partial summation at most
\begin{equation}\label{case2}
\ll X^{\varepsilon} \sum_{\substack{z^{1/2} \leq a \leq z \\ P^+(a) \leq \xi\\ (a, q) = 1}} S(a, 1, 1, 1) \ll \frac{X^{1+\varepsilon}}{ z^{1/2} } + X^{\varepsilon} Y.
\end{equation}
The same argument works with exchanged indices if $n_2$ is in case (III). 

3) Finally we consider the situation when both $n_1$ and $n_2$ are in cases (I) or (IV). 

Suppose first that $n_1$ is in case (I). Then $  z^{\Omega(b_1)/2} \leq b_1 \ll X^{O(1)}$, so that $\Omega(b_1) \ll 1$ and hence $\lambda_1(b_1) \ll 1$.  We conclude that  $n_1$ is divisible by some $$a_1 \leq  z, \quad (n_1/a_1, \mathcal{P}_{q, z^{1/2}}) = 1$$ and $\lambda_1(n_1) \ll \lambda_1(a_1)$. The same argument works if $n_2$ is in case (I) with exchanged indices. 

Next suppose that $n_1$ is in case (IV). Then we localize $p_{k_1+1}$ into intervals of the shape $z^{1/(r+1)} \leq p_{k_1+1} \leq z^{1/r}$ for $2 \leq r \leq \log z/\log \xi$. In each such interval we  have $p_{k_1} \leq p_{k_1 + 1} \leq z^{1/r}$. We conclude that for each fixed $r$, the number $n_1$ is divisible by some $$z^{1/2} < a_1 \leq z, \quad P^+(a_1) \leq z^{1/r}, \quad (n_1/a_1, P_{q, z^{1/(r+1)}}) = 1,$$
and since $\Omega(n_1/a_1) \leq r$, we have $\lambda_1(n_1) \ll \lambda_1(a_1) k^r$ since $\lambda_j(n) \leq \tau_k(n)$. The same argument works if $n_2$ is in case (IV) with exchanged indices. 

We conclude that in either case the contribution of such pairs is
\begin{displaymath}
\begin{split}
& \sum_{ r_1, r_2 \leq \log z/\log \xi} k^{r_1 + r_2} \sum_{\substack{\delta_{r_1 \not= 1} z^{1/2} \leq a_1  \leq z\\ \delta_{r_2 \not= 1} z^{1/2} \leq a_2  \leq z\\ P^+(a_1), P^+(a_2) \leq z^{1/r} \\ (a_1a_2, q) = 1}} \lambda_1(a_1) \lambda_2(a_2) S(a_1, a_2, z^{1/(r_1 + 1)}, z^{1/(r_2 + 1)})\\
& \ll  \sum_{ r_1, r_2 \leq \log z/\log \xi} k^{r_1 + r_2}r_1r_2 \sum_{\substack{\delta_{r_1 \not= 1} z^{1/2} \leq a_1  \leq z\\ \delta_{r_2 \not= 1} z^{1/2} \leq a_2  \leq z\\ P^+(a_1), P^+(a_2) \leq z^{1/r}}} \frac{X\lambda_1(a_1) \lambda_2(a_2)  }{  a_1a_2 (\log z)^2} + X^{\varepsilon}Yz^3 
\end{split}
\end{displaymath}
where we used once again Lemma \ref{lem3}. By \cite{Wo}*{Lemma 3} we have
$$ \sum_{\substack{\delta_{r_1 \not= 1} z^{1/2} \leq a   \leq z \\ P^+(a )  \leq z^{1/r}}} \frac{ \lambda_j(a)   }{a} \ll \exp\Big(\sum_{p \leq z} \frac{\lambda_j(p)}{p} - c r \log r\Big)$$
for some constant $c > 0$, $r \ll \log z/\log\log z$ (this condition is needed in the proof and requires us to treat case (III) separately) and $j = 1, 2$, 
 so that the total contribution of pairs in the present case  is
\begin{equation}\label{case3}
\ll \frac{X}{  (\log z)^2} 
 \exp\Big(\sum_{p \leq z} \frac{\lambda_1(p) + \lambda_2(p)}{p}\Big)  
+ X^{\varepsilon}Yz^3.
\end{equation}

Combining \eqref{case1}, \eqref{case2} and \eqref{case3}, we complete the proof. \qed\\

We will apply Proposition \ref{thm-sieve} to the multiplicative functions $n\mapsto |\lambda_{f_i}(n)|$, $i=1,2$, 
where $f_1, f_2$ are cuspidal (holomorphic or Maa{\ss}) Hecke eigenforms with Hecke eigenvalues $\lambda_{f_i}(n)$  
satisfying the Ramanujan-Petersson conjecture, so that 
\begin{equation}
	\label{RPbound}
|\lambda_{f_1}(n)|,|\lambda_{f_2}(n)|\leq \tau(n)
\end{equation}
for all $n \geq 1$. 

\begin{lemma}\label{sato} Let $f$ be a  cuspidal newform with Hecke eigenvalues $\lambda_f$  satisfying the Ramanujan-Petersson conjecture.

-- If $f$ is a dihedral form (in which case the Ramanujan-Petersson conjecture is automatic), then for any $z\geq 5$ 
we have 
\begin{equation}
	\label{CMbound} \sum_{p \leq z} \frac{|\lambda_f(p)|}{p}  \leq  
 \frac{3}{4} \log\log z + O\big(1 + \log (\text{{\rm cond}} f)\big).
	\end{equation}

-- If $f$ is not a dihedral form, then for any $z\geq 5$ we have
\begin{equation}
	\label{nonCMbound}
	\sum_{p \leq z} \frac{|\lambda_f(p)|}{p}  \leq  \frac{17}{18} \log\log z + O\big(1 + \log (\text{{\rm cond}} f)\big).
\end{equation}
\end{lemma}

\textbf{Proof.} Let $N$ denote the level of $f$ and $\chi$ be its nebentypus.  
If $f$ is dihedral, we recall that it  comes from some Hecke character (not of order $1$ or $2$ since $f$ is cuspidal) of a quadratic number field $K/\mathbb{Q}$ and hence $\lambda_f(p)$ vanishes if $p$ is inert in $K$. For prime $p\nmid N$ split in $K$ we have
$$|\lambda_f(p)|\leq \frac{1+|\lambda_f(p)|^2}{2}=\frac32+\frac{1}{2}\lambda_{f_2}(p)$$
where $f_2$ is the  theta series induced from the square of the Hecke character mentioned above.
We conclude that
 $$ \sum_{p \leq z} \frac{|\lambda_f(p)|}{p}  \leq  \frac32\sum_{\substack{p \leq z,(p,N)=1\\ p \text{ split}}} \frac{1}{p}w\Big(\frac{p}{z}\Big) + \frac{1}{2}\sum_{\substack{p \leq z,(p,N)=1\\ p \text{ split}}} \frac{\lambda_{f_2}(p)}{p}w\Big(\frac{p}{z}\Big)+O(\log N)$$
 for  a smooth non-negative function $w$ supported on $[0,2)$ and equal to $1$ on $[0,1]$.
 
 The desired bound \eqref{CMbound} follows after integration by parts from Perron's formula applied to $-\mathrm{d}\log \zeta_K(s)$ and $-\mathrm{d}\log L(s,f_2)$ which are holomorphic and logarithmically bounded in a standard (Hadamard/de la Vall\'ee-Poussin type) neighbourhood of the line $\Re s=1$ except   for a pole at $s=1$ (for the former) with residue $+1$  and (possibly) a pole at some Siegel zero with residue $-1$ (which if it exists, then contributes   a negative amount).

	If $f$ is not dihedral, the bound \eqref{nonCMbound} is similar, but requires a few higher symmetric power $L$-functions of $f$  (cf.\ \cite{EMS}). The starting point is the following computation, valid for $p\nmid N$ and whenever $|\lambda_f(p)| \leq 2$, namely
$$|\lambda_f(p)| \leq 1 + \frac{1}{2}\lambda_f(p^2)\ov\chi(p) - \frac{1}{18} \lambda_f(p^2)^2\ov\chi^2(p) = \frac{17}{18} + \frac{4}{9} \lambda_{\sym^2 f}(p)\ov\chi(p) - \frac{1}{18} \lambda_{\sym^4 f}(p)\ov\chi^2(p).$$
We have therefore
\begin{displaymath}
\begin{split}
 \sum_{\substack{p \leq z\\ p \nmid N}} \frac{|\lambda_f(p)|}{p} & \leq  \sum_{p \nmid N  } \frac{|\lambda_f(p)|}{p} w\Big(\frac{p}{z}\Big)\\
 & \leq \frac{17}{18} \log\log z + O(1) + \frac{4}{9}  \sum_{p  \nmid N } \frac{\lambda_{\sym^2f}(p)\ov\chi(p)}{p} w\Big(\frac{p}{z}\Big)  - \frac{1}{18}   \sum_{p \nmid N } \frac{\lambda_{\sym^4f}(p)\ov\chi^2(p)}{p} w\Big(\frac{p}{z}\Big).
 \end{split}
 \end{displaymath}
The bounds for the last two sums are then consequences of the holomorphy and the logarithmic bounds satisfied by $-\mathrm{d}\log L(s, \sym^{2\nu} f\times \ov\chi^\nu)$ for $\nu=1,2$ when $s$ is within the standard neighbourhood of the line $\Re s=1$.
 
For $\nu=1$, the required properties follow from the cuspidality of $L(s,\sym^2 f\times \ov\chi)$ (since $f$ is non-dihedral),  and the absence of a Siegel zero established in \cites{HL,Ba}. For $\nu=2$, we know from \cites{KimSar,Kimsh}  that $L(s,\sym^4 f\times \ov\chi^2)$ is automorphic but not necessarily cuspidal. It is cuspidal unless $f$ is tetrahedral or octahedral. In the tetrahedral case, we have by \cite{Kimsh}*{\S 3.2}
$$L(s, \sym^4 f\times\ov\chi^2)=L(s, \sym^2 f\times \ov\chi)L(s,\chi_3)L(s,\ov\chi_3)$$
for $\chi_3=\chi_{3,f}$ a non-trivial character of order $3$. In particular $L(s, \sym^4 f\times \ov\chi^2)$ has no zeros or poles in a standard region along the line $\Re s=1$. In the octahedral case we have
\cite{Kimsh}*{Thm.\ 3.3.7}
$$L(s, \sym^4 f\times \ov\chi^2)=L(s,\pi(\chi))L(s,\sym^2 f\times \eta\chi_f^{-1})$$
 where $\eta=\eta_f$ is a quadratic character  and $\pi(\chi)$ is the dihedral representation induced from a non-trivial cubic Hecke character $\chi=\chi_f$ of the quadratic field determined by $\eta$. It follows that $L(s,\sym^4 f\times \ov\chi^2)$ has no zeros or poles in a standard region along the line $\Re s=1$.

 Finally if $f$ is not of the above type then $L(s, \sym^4 f\times \ov\chi^2)$ is cuspidal and has no zeros or poles in a standard region along the line $\Re s=1$ except for a possible Siegel zero. If $\chi_f$ is trivial, the existence of a Siegel zero was ruled out by Ramakrishnan and Wang in \cite{RW}*{Theorem B'}. As was pointed out to us by D.\  Ramakrishnan,  the proof extends to the case of a general nebentypus with adequate modifications: let $\Pi$ be the degree $9$ isobaric sum 
 $$\Pi=1\boxplus \sym^2f\times\ov\chi\boxplus \sym^4f\times\ov\chi^2;$$
 its degree $81$ Rankin-Selberg $L$-function (which has non-negative coefficients) factors as 
\begin{align}\nonumber
    L(s,\Pi\times \Pi)=&\zeta(s)L(s,\sym^2f\times\ov\chi)^4L(s,\sym^4f\times\ov\chi^2)^4 L(s,\sym^6f\times\ov\chi^3)^2\\
   \label{LPi} &
 \times L(s,\sym^2f\times\sym^2f\times\ov\chi^2)L(s,\sym^4f\times\sym^4f\times\ov\chi^4)\\
 =&L(s,\sym^4f\times\ov\chi^2)^4L_2(s)\nonumber
\end{align}
say. Following the proof of \cite{RW}*{Theorem B'} it is sufficient to prove that $L_2(s)$ is holomorphic along the interval $(1/2,1)$ (cf. \cite{RW}*{Prop. 5.21}). Using the factorisation
$$L(s,\sym^3f;\sym^2\times\ov\chi^3)=L(s,\sym^2f\otimes\ov\chi)L(s,\sym^6f\otimes\ov\chi^3),$$
one can rewrite $L_2(s)$ into the form
\begin{multline*} L_2(s)=\zeta(s)L(s,\sym^2f\times\ov\chi)^2
L(s,\sym^3f;\sym^2\times\ov\chi^3)^2\\
\times L(s,\sym^2f\times\sym^2f\times\ov\chi^2)L(s,\sym^4f\times\sym^4f\times\ov\chi^4), 
\end{multline*}
 and its remains to prove the holomorphy of $L(s,\sym^3f;\sym^2\times\ov\chi^3)$ along $(1/2,1)$. 
One then proceeds as in \cite{RW}*{Lem. 5.25 \& \S 7}. The only difference is that, when $\chi^3$ is non-trivial, one has to use the work of Takeda \cite{Tak} in place of the work of Bump-Ginzburg \cite{BG} regarding the holomorphy of the incomplete $L$-function $L^S(s,\sym^3f;\sym^2\times\ov\chi^3)$ (where $S$ is the union of the archimedean places and the finite places where $f$ is ramified). 
 \qed
 
 \begin{remark}
     In this paper we will need only the case where $f$ has trivial nebentypus.
 \end{remark}
 
 \begin{remark} In the non-CM case the Sato-Tate conjecture predicts the constant $8/(3\pi)$ in place of 17/18 in \eqref{nonCMbound}, and it is known when $f$ is holomorphic \cites{CHT,HST,Ta,NewTh}. In the CM case, the best possible constant is $2/3$ which is attained if $f$ comes from a character of order 3. 
  \end{remark}

\begin{cor} Let $f_1,\ f_2$ be cuspidal    newforms with Hecke eigenvalues $(\lambda_{f_i}(n))_{n\geq 1}$, $i=1,2$, satisfying the Ramanujan-Petersson conjecture \eqref{RPbound}. Then, with the notation  of Proposition 	\ref{thm-sieve} we have
\begin{equation*}
	\sum_{\substack{ (n_1, n_2) \in \mathcal{S}\\ (n_1n_2, q) = 1}} |\lambda_{f_1}(n_1) \lambda_{f_2}(n_2)|  \ll_{f_1,f_2 ,\gamma,\varepsilon} \frac{X}{(\log z)^{1/9}}
+  \frac{X^{1+\varepsilon}}{ z^{1/4}} + X^{\varepsilon} Y z^3
\end{equation*}
for any $\eps>0$.
\end{cor}

\section{Equidistribution for the cuspidal spectrum}

In this section, we establish the equidistribution statements of Theorems \ref{main1} and \ref{main2} for pairs of cuspidal newforms (assuming the Ramanujan-Petersson conjecture).

\subsection{The discrete case}

For $q\in \mathbb{N}$  and $b \in (\mathbb{Z}/q\mathbb{Z})^{\times}$ we recall the definition \eqref{lattice} of the lattice
$$\Lambda_{q;  b}=\{(n_1,n_2)\in\Zz^2,\ 
n_1 + b n_2 \equiv 0 \mods q\}$$ 
and its minimum  $s = s(b; q)$. We start with the following simple result.

\begin{lemma}\label{min} Let $(db, q) = 1$. Then 
$$d^{-1} s(b; q) \leq  s(bd; q) \leq d s(b; q) $$
and
$$d^{-1} s(b; q) \leq  s(b\bar{d}; q) \leq d s(b; q) $$
and $s(b; q) = s(\bar{b}; q)$. 
    \end{lemma}

\textbf{Proof.} Clearly $(n_1, n_2) \in \Lambda_{q; b}$ implies $(n_1 d, n_2) \in \Lambda_{q; bd}$, so that $s(bd; q) \leq d s(b; q)$. Similarly, $s(b\bar{d}; q) \leq d s(b; q)$. Replacing $b$ with $b\bar{d}$ and $bd$ we obtain the inequalities in the other direction. The last statement follows from exchanging $n_1$ and $n_2$. \qed\\

With future applications in mind, we consider slightly more general Weyl sums than in \eqref{Weylsum0}.

\begin{theorem}\label{thm0} Let $q$ be a large prime, $b \in (\mathbb{Z}/q\mathbb{Z})^{\times}$ and $s = s(b; q)$. 

Let $f_j$, $j=1,2$, be two $L^2$-normalized cuspidal  Hecke-Maa{\ss} newforms of  levels $N_j$ with Hecke eigenvalues $\lambda_j(n)$ for which the Ramanujan-Petersson conjecture holds. 

Let $x_0 \in \mathbb{R}$, $y_0 \in \mathbb{R}^{\times}$, $r_0 \in \mathbb{Q}$ where the denominator of $r_0$ is coprime to $N_1N_2q$, and let
\begin{equation}\label{defW} \mathcal{W}_{f_1, f_2}(b; q; x_0, y_0, r_0) :=  \frac{1}{q}\sum_{a\mods q} f_1\Big( \frac{a + i}{q}\Big) f_2\Big( \frac{ba + x_0 + y_0i}{q} + r_0\Big).
\end{equation}
One has
$$\mathcal{W}_{f_1, f_2}(b; q; x_0, y_0, r_0)\ll_{f_1, f_2, x_0, y_0, r_0, \varepsilon} s^{\varepsilon - 1}
 + (\log q)^{-1/9}$$
for any $\varepsilon > 0$ with an implied constant depending polynomially on $|x_0|$, $y_0 + y_0^{-1}$, the denominator of $r_0$ and the conductors of $f_1$ and $f_2$.  In other words, for fixed $(x_0, y_0, r_0)$ and fixed cusp forms $f_1, f_2$ we obtain decay  as soon as $s \rightarrow \infty$. 
\end{theorem}

\begin{remark}
The reason to include general triples $(x_0, y_0, r_0)$ (and not just the triple $(0, 1, 0)$ as in \eqref{Weylsum0}) is that we will later apply this to pairs of cusp forms $(f_1, f_2)$ where $f_2$, lifted to a function on an $S$-adic quotient of $${\rm PGL}_2(\mathbb{R}) \times {\rm PGL}_2(\mathbb{Q}_{q_1}) \times {\rm PGL}_2(\mathbb{Q}_{q_2})$$ for two fixed primes $q_1, q_2$, is acted on by a fixed group element $g \in  {\rm PGL}_2(\mathbb{R}) \times {\rm PGL}_2(\mathbb{Q}_{q_1}) \times {\rm PGL}_2(\mathbb{Q}_{q_2})$. When translated back to classical language, this action has the effect of introducing the extra parameters $(x_0, y_0, r_0)$. See Section \ref{sec5} for details.

We remark that a similar bound can be obtained for more general test functions, for instance Maa{\ss} forms of fixed weights $k_j$, or even more generally automorphic forms whose archimedean component is a fixed test function (for instance compactly supported) in the Kirillov model.  
\end{remark}

\begin{remark}\label{sisoptimal} 
 To see  that the condition $s \rightarrow \infty$ cannot be dropped completely or replaced with a simpler condition of the kind $q \| b/q \| \rightarrow \infty$,   consider the case where $b = (q+1)/2$ (for odd $q$). Then the Fourier expansion \eqref{sast} below yields essentially
$$\frac{1}{q} \sum_{\substack{n_1, n_2 \asymp q\\ q \mid n_1 - bn_2}} \lambda_1(n_1) \lambda_2(n_2).$$
The congruence is  equivalent to $q \mid 2n_1 - n_2$, and so we obtain a diagonal term
$$\frac{1}{q} \sum_{ n \asymp q} \lambda_1(a) \lambda_2(2a)$$ 
which in the case $f_1 = f_2$ and $\lambda_2(2) \not= 0$ does not decay in $q$. In other words, the underlying condition for equidistribution is really of diophantine nature.  
\end{remark}

\textbf{Proof.} We start with the Fourier expansion
\begin{equation}
	\label{Fourierexpansion}
	f_j(z) =   \sqrt{y}  \sum_{n \geq 1}  \lambda_j(n) \Big( \frac{2 \cosh(\pi t_j)}{L(1, \sym^2 f_j)}\Big)^{1/2} 
   K_{it_j}(2\pi n y) \big(e(nx)+ \epsilon_j e(-nx)\big)
\end{equation}	 
where $\epsilon_j$ is the parity, $\lambda_j(n)$ is the $n$-th Hecke eigenvalue (we have assumed that the $f_j$ are newforms)  and\footnote{by the Selberg eigenvalue conjecture, although this is not essential for the argument, unlike the Ramanujan-Petersson conjecture at finite places}  
$t_j \geq 0$ is the spectral parameter of $f_j$ for $j = 1, 2$.  
For notational simplicity let us write $$L_j = L(1, \sym^2f_j)^{1/2}, \quad K^{\ast}_{it}(x) = \cosh(\pi t)^{1/2} K_{it}(2 \pi x).$$
We  use the simple bound 
$$x^j \frac{d^j}{dx^j} K^{\ast}_{it}(x) \ll_{j, A, \varepsilon, t} x^{- \varepsilon}(1 + x)^{-A} $$ 
for $j\in \mathbb{N}_0$, $A, \varepsilon  > 0$ with 
polynomial dependence in $t$. In the following we use the convention that all implied constants may depend polynomially on $\text{cond}(f_j)$ and $|x_0|$, $y_0 + y_0^{-1}$, $\text{den}(r_0)$ without further mention.  

Summing over $a$ (mod $q$), we conclude 
\begin{displaymath}
\begin{split}
\mathcal{W}_{f_1, f_2}(b; q; x_0, y_0, r_0) & = \frac{2|y_0|^{1/2}}{L_1L_2 q}  
\Big( (1 + \epsilon_1\epsilon_2)\sum_{q \mid n_1 + bn_2} +(\epsilon_1+\epsilon_2)  \sum_{q \mid n_1 -bn_2} \Big)\lambda_1(n_1) \lambda_2(n_2)\\
& \times e\Big(\frac{n_2x_0}{q} + n_2r_0\Big) K^{\ast}_{it_1}\Big( \frac{n_1}{q} \Big)K^{\ast}_{it_2}\Big(  \frac{|y_0|n_2}{q}\Big) . 
\end{split}
\end{displaymath}
This vanishes unless $\epsilon_1=\epsilon_2 = \epsilon$, say,  in which case we get
\begin{equation*}
\begin{split}
\mathcal{W}_{f_1, f_2}(b; q; x_0, y_0, r_0) &= 
 \frac{1 }{ q}  \sum_{\pm}  \sum_{ q \mid n_1 \pm  bn_2  }  \lambda_1(n_1) \lambda_2(n_2) e(n_2r_0)G\Big(\frac{n_1}{q}, \frac{n_2}{q}\Big)  
\end{split}
\end{equation*}
 with
\begin{equation}\label{G}
G(x_1, x_2) =\frac{4e(x_0x_2)}{L_1L_2   }  K^{\ast}_{it_1}( x_1)K^{\ast}_{it_2}(|y_0|x_2) 
\end{equation}
satisfying $$|x_1|^{j_1} |x_2|^{j_2} \frac{d^{j_1}}{dx_1^{j_1}}\frac{d^{j_2}}{dx_2^{j_2}} G(x_1, x_2) \ll_{\varepsilon, A, j_1, j_2,x_0,y_0} |x_1x_2|^{-\varepsilon}(1 + |x_1| + |x_2|)^{-A}.$$
By trivial estimates we can truncate the $n_1, n_2$-sum at $q^{1+\varepsilon}$ at the cost of a negligible error. 
Moreover, if $r_0 = c_0/d_0$, we can replace the additive character $n_2 \mapsto e(n_2r_0)$ with (a linear combination of) Dirichlet characters $\chi$ modulo $d_0$,  
so that it suffices to bound 
 \begin{equation}\label{sast}
\begin{split}
 \frac{1 }{ q}    \Big|\sum_{\substack{q \mid n_1 \pm  bn_2  \\ n_1, n_2 \ll q^{1+\varepsilon}}}  \lambda_1(n_1) \lambda_2(n_2)\chi(n_2) G\Big(\frac{n_1}{q}, \frac{n_2}{q}\Big)\Big| 
\end{split}
\end{equation}
for some  choice of $\pm$. 

The character $\chi$ may not be primitive; if it is induced by the primitive character $\chi^{\ast}$, we obtain by M\"obius inversion and the Hecke multiplicativity relations 
\begin{equation*} 
\begin{split}
 &\frac{1 }{ q}    \Big|\sum_{\substack{q \mid n_1 \pm  bn_2  \\ n_1, n_2 \ll q^{1+\varepsilon}}}  \lambda_1(n_1) \lambda_2(n_2)\chi^{\ast}(n_2)\sum_{f \mid (n_2, d_0)} \mu(f) G\Big(\frac{n_1}{q}, \frac{n_2}{q}\Big)\Big|   
\\ \leq  &\frac{1 }{ q} \sum_{f \mid d_0}   \Big|\sum_{\substack{q \mid n_1 \pm  bfn_2  \\ n_1, fn_2 \ll q^{1+\varepsilon}}}  \lambda_1(n_1) \sum_{g \mid (f, n_2)} \mu(g)\lambda_2\Big(\frac{f}{g}\Big) \lambda_2\Big(\frac{n_2}{g}\Big)\chi^{\ast}(n_2)  G\Big(\frac{n_1}{q}, \frac{fn_2}{q}\Big)\Big| \\
\\ \leq &\frac{1 }{ q} \sum_{g \mid f \mid d_0} \Big|\lambda_2\Big( \frac{f}{g}\Big)\Big|   \Big|\sum_{\substack{q \mid n_1 \pm  bfgn_2  \\ n_1, n_2fg \ll q^{1+\varepsilon}}}  \lambda_1(n_1)  \lambda_2(n_2)\chi^{\ast}(n_2)  G\Big(\frac{n_1}{q}, \frac{gfn_2}{q}\Big)\Big| .
\end{split}
\end{equation*}
The point of this maneuvre is to show that in \eqref{sast} we may assume without loss of generality that $\chi$ is primitive. Since $(d_0, N_2) = 1$, the function $n \mapsto \lambda_2\chi(n)$ describes the Hecke eigenvalues of the \emph{newform} $f_2 \times \chi$, which allows us to apply Proposition \ref{thm1}. To this end we have to replace $b$ by some fixed multiple, but by Lemma \ref{min} we have $s(b; q) \asymp s(fgb; q)$, so that the subsequent analysis remains unchanged. We may therefore return to \eqref{sast} under the assumption that $\chi$ is primitive.

The $(n_1, n_2)$-sum runs through a lattice $\Lambda_{q,\pm b}$ with basis
$$\binom{\mp b}{1}, \quad \binom{q}{0}$$
of volume $q$ and minima $0<s(b;q)=s(-b;q) \ll q^{1/2}$.  Let us fix one sign and drop it from the notation. 
Let 
$$\binom{x_1}{x_2}, \quad \binom{y_1}{y_2},$$
  be a reduced basis of $\Lambda_{q, b}$ with the first vector of minimal length so that $|x_1| + |x_2| \ll s$. Both $x_1, x_2$ are nonzero (since $q \nmid b$) and coprime (since $s$ is minimal). In terms of this basis, the congruence condition reads
\begin{equation}\label{cong}
n_1x_2 \equiv n_2x_1 \, (\text{mod } q).
\end{equation}

We now use two different methods to estimate the $n_1, n_2$-sum depending on the size of $s$.

 If $s$ is not too big, then we interpret the $n_1, n_2$-sum as a shifted convolution problem and apply  first Proposition \ref{thm1} to the inner sum with summation condition \eqref{cong} which we write as an equality
 $$n_1x_2 - n_2x_1 = h$$ with $q \mid h$ (see \cite{BFKMM} for a similar argument). 
 
 Since $s \ll q^{1/2}$ and $n_1, n_2 \ll q^{1+\varepsilon}$, we have automatically that  $(h/q, q) = 1$. We need to treat the diagonal term $\Delta$ with $h = 0$ separately. In this case $x_1 \mid n_1$, $x_2 \mid n_2$, and by Rankin-Selberg theory (and the Ramanujan conjecture, although this could be dispensed with) it is easy to see the contribution is
\begin{equation}\label{diag}
 \Delta  \ll \frac{1}{s^{1-\varepsilon}}.
\end{equation}

Let  us call $\mathcal{W}^{\ast}_{f_1, f_2}(b; q;x_0,y_0, r_0) $ the remaining portion of \eqref{sast}. 
Applying a smooth partition of unity to the $n_1, n_2, h$-sums, we obtain by Proposition \ref{thm1} with $$d = q, \quad M_1, M_2 \ll q^{1+\varepsilon}, \quad H \ll l_1M_1+ l_2M_2,  \quad l_1, l_2 \ll s$$  the upper bound
\begin{equation}\label{sast1}
\mathcal{W}^{\ast}_{f_1, f_2}(b; q; x_0, y_0, r_0) \ll
s^{1+\theta} q^{\theta - 1/2 + \varepsilon}.
\end{equation}
Note that the $\theta$-dependence in Proposition \ref{thm1} comes from the entire spectrum of forms of level $x_1x_2$, not from the two forms $f_1, f_2$. 

Alternatively, we apply absolute values to the inner sum in \eqref{sast}, and we can then apply Proposition \ref{thm-sieve} with $z = q^\gamma$ for some very small $\gamma > 0$ and $$\mathcal{S} = \Lambda \cap \mathcal{R}$$ for suitable sets $\mathcal{R} \subseteq \mathbb{R}_{>0}^2$. More precisely, we estimate the contribution of the terms with $q \mid n_1n_2$ trivially by $O(q^{\varepsilon - 1})$. For the remaining terms  we apply a decomposition into annuli $(k-1)q \leq \sqrt{n_1^2 + n_2^2} \leq k q$ for $k = 1, 2, \ldots, q^{\varepsilon}$, and by Proposition \ref{thm-sieve} we obtain the bound
$$\mathcal{W}_{f_1, f_2}(b; q; x_0, y_0, r_0) \ll   q^{\varepsilon} \Big(\frac{1}{q}   + \frac{1}{z^{1/4}} + \frac{z^3}{s}\Big) +   \frac{1}{(\log q)^2} \exp\Big( \sum_{p \leq z} \frac{|\lambda_1(p)| + |\lambda_2(p)|}{p} \Big).$$
Employing Lemma \ref{sato}, we obtain
\begin{equation}\label{sast2}
\mathcal{W}_{f_1, f_2}(b; q; x_0, y_0, r_0) \ll    \frac{1}{(\log q)^{1/9}}   + \frac{q^{3\gamma + \varepsilon}}{s} .
\end{equation}
Combining \eqref{diag}, \eqref{sast1} and \eqref{sast2}, we complete the proof of Theorem \ref{thm0}. \qed

\begin{remark} The bound \eqref{sast2} is almost sufficient alone. The shifted convolution argument is only needed when $s$ grows below the $q^{\varepsilon}$-scale, say logarithmically with $q$.  On the other hand, however, the shifted convolution argument is more robust and provides a power saving.
\end{remark}

\subsection{The continuous case}

For the continuous joint equidistribution problem in Theorem \ref{main2} we need the following analogue of Theorem \ref{thm0}. 

\begin{theorem}\label{thm-cont}  Let $f_1, f_2$ be two $L^2$-normalized cuspidal Hecke-Maa{\ss}  newforms of   levels $N_j$ with Hecke eigenvalues $\lambda_j(n)$ for which the Ramanujan-Petersson conjecture holds. Let $x_0 \in \mathbb{R}$, $y_0 \in \mathbb{R}^{\times}, r_0 \in \mathbb{Q}$ such that the denominator of $r_0$ is coprime to $N_1N_2$.  Let $T> 1$, $I \subseteq (0, \infty)$ a fixed compact interval, $y \in I$ and $W$ a fixed smooth weight function with compact support in $\mathbb{R}$ and write $y = a/q + O(1/qQ)$ for positive coprime integers $a, q$ with $q \leq Q := T^{0.99}$. Then
$$\int W(x) f_1\Big(x + \frac{i}{T}\Big) f_2\Big(xy  + r_0+ \frac{x_0 + iy_0}{T}\Big) dx \ll_{\varepsilon, f_1, f_2, x_0, y_0, r_0} \frac{1}{(\log T)^{1/9}} + \frac{1}{q^{\varepsilon - 1}}$$
for any $\varepsilon > 0$ with an implied constant depending polynomially on the conductors of $f_1$ and $f_2$ and on $|x_0|$, $y_0 + 1/y_0$ and the denominator of $r_0$. 
\end{theorem}

\textbf{Proof.} We start again with the Fourier expansion \eqref{Fourierexpansion} and perform the $x$-integration.  As in \eqref{sast} this leaves us with bounding
$$\mathcal{R}_{f_1, f_2}(y, T) :=  \frac{1}{T} \sum_{\pm} \Big| \sum_{ n_1, n_2 \ll T^{1+\varepsilon} } \hat{W}(n_1 \pm n_2y)\lambda_1(n_1) \lambda_2(n_2) \chi(n_2) G\Big(\frac{n_1}{T}, \frac{n_2}{T}\Big)  \Big| $$
where $\hat{W}$ denotes the Fourier transform of $W$, $G$ is as in \eqref{G} and $\chi$ is some primitive character whose conductor divides the denominator of $r_0$.   Inserting the rational approximation, the $n_1, n_2$-sum is, up to a negligible error, restricted by $$qn_1 \pm an_2 \ll H:= T^{\varepsilon} \Big(q + \frac{T}{Q}\Big). $$
For now, let us only assume  $\log Q \asymp \log T$. Applying Proposition \ref{thm1} with $d = 1$ and treating the diagonal contribution $qn_1 \pm an_2 = 0$ separately as in \eqref{diag}, we obtain the first bound
\begin{equation}\label{R1}
\begin{split}
\mathcal{R}_{f_1, f_2}(y, T) & \ll   \frac{1}{q^{1-\varepsilon}} + \frac{T^{\varepsilon}}{T}(qT)^{1/2 + \theta} \Big(q + \frac{T}{Q}\Big)^{1/2}\Big( 1+ \frac{T}{qQ}\Big)^{1/2}\\
&\ll \frac{1}{q^{1-\varepsilon}} + T^{\varepsilon}(qT)^{\theta}\Big( \frac{q}{T^{1/2}} + \frac{T^{1/2}}{Q}\Big).
\end{split}
\end{equation}
 On the other hand, we can apply absolute values to $\mathcal{R}_{f_1, f_2}(y, T)$ and use Proposition \ref{thm-sieve}. To this end,  let $\mathcal{S}_{\pm}(a, q, H', T')$ for $H' \in \mathbb{N}$ and $T' \geq 1$ be the number of $(n_1, n_2)$ such that $T' < n_2 \leq 2T' $  and $H' < qn_1 \pm a n_2 \leq 2H'$. Then
 $$\sum_{\substack{\lambda \in \mathcal{S}_{\pm}(a, q, H', T') \\ d_1 \mid \lambda_1, d_2 \mid \lambda_2}} 1 = H'\Big(\frac{T'}{qd_1d_2} + O(1)\Big),$$
so that
$$\mathcal{S}_{\pm}(a, q, H', T')\in \mathcal{C}_1(H'T'/q, H').$$
Recall that $\widehat{W}(x) G(x_1, x_2)\ll_A (1 + x)^{-A} |x_1x_2|^{-\varepsilon}(1 + |x_1|)^{-A} (1 + |x_2|)^{-A}$.  
Now from Theorem \ref{thm-sieve} and Lemma \ref{sato} we conclude for some sufficiently small $\gamma > 0$ similarly as in the preceding proof
\begin{equation}\label{R2}
\begin{split}
\mathcal{R}_{f_1, f_2}(y, T) & \ll T^{-10} + \frac{1}{T} \Bigg( \frac{(q + T/Q)T}{q (\log T)^{1/9}}\Big) + T^{\varepsilon}\Big(\frac{(q + T/Q)T}{q T^{\gamma/4}}\Big) + \Big(q + \frac{T}{Q}\Big)T^{3\gamma}\Bigg)\\
&\ll \frac{1}{(\log T)^{1/9}} + \frac{T}{qQ} + \frac{q}{T^{1 - 3\gamma - \varepsilon}} .
\end{split}
\end{equation}
We may now choose, for instance, $Q = T^{0.99}$ (with $\gamma$ as above sufficiently small) and combine \eqref{R1} and \eqref{R2} to obtain the desired bound
$$\mathcal{R}_{f_1, f_2}(y, T) \ll \frac{1}{(\log T)^{1/9}} + \frac{1}{q^{\varepsilon - 1}} + T^{\varepsilon} \min_{q \leq Q} \Big((qT)^{\theta} \frac{q}{T}, \frac{T^{4\gamma}}{q}\Big),$$
and the last term is absorbed in the previous two terms.

\section{Application of the joinings theorem of Einsiedler-Lindenstrauss}\label{sec5}

In this section we use the powerful measure classification theorem of Einsiedler-Lindenstrauss \cite{EL}.

\subsection{Passage to an $S$-arithmetic quotient}

For any place $v$  of $\Qq$ we denote by $\Qq_v$ the corresponding completion and for a finite place $v=p$ by $\Zz_p\subseteq \Qp$ the closure of $\Zz$. We denote by $$\Aa=\Rr\times \Aa_f=\Rr\times{\prod_{p}}'\Qp$$
the ring of adeles (the first component is always the infinite place). 

We fix an integer $D\geq 1$ and a non-empty finite set of places $S=\{\infty\}\cup S_f$  containing $\infty$ and such that the primes in $S$ are coprime to $D$; we set
$$\Qq_S=\prod_{v\in S}\Qv,\quad \Zz_S=\prod_{p\in S}\Zz_p,\quad \Zz[1/S]=\Zz\Big[\prod_{p\in S}\frac{1}p\Big]
,$$
and
$$\Aa^{(S)}={\prod_{v\not\in S}}'\Qv,\quad \whZ^{(S)}=\prod_{p\not\in S}\Zp.$$

Let $\rmG=\PGL_2$ and let $K(D)^{(S)} $ be the open-compact subgroup
$$K(D)^{(S)}=\prod_{p\not\in S}K(D)_p\subseteq \rmG(\Aa^{(S)})$$
where for $p\not\in S$ we put 
$$K(D)_p=\Big\{\begin{pmatrix}
    a&b\\c&d
\end{pmatrix}\in \rmG(\Zp),\ c\in D\Zp\Big\}\subseteq\rmG(\Zp); $$
 in particular one has $K(D)_p=\rmG(\Zp)$ for $p\nmid D$.
 We also set
 $$K_\infty=\PSO_2(\Rr)\subset \rmG(\Rr).$$
 
 By the strong approximation property for $\SL_2$ and the fact that $\det(K(D)^{(S)})=\whZ^{(S)\times}/(\whZ^{(S)\times})^2$ we have 
$$\rmG(\Qq)\rmG(\Qq_S)K(D)^{(S)}=\rmG(\Aa)$$ and
$$\rmG(\Qq)\bash \rmG(\Aa)/K(D)^{(S)}\simeq \Gamma(D)_S\bash\rmG(\Qq_S),$$
where
$$\Gamma(D)_S:=\rmG(\Qq)\cap K(D)^{(S)}=\Big\{\begin{pmatrix}
    a&b\\c&d
\end{pmatrix}\in \rmG(\Zz[1/S]),\ c\in D\Zz[1/S]\Big\}\subseteq \rmG(\Zz[1/S])$$ is a discrete subgroup of $\rmG(\Qq_S)$. In the sequel we set
$$X(D)^S=\Gamma(D)_S\bash\rmG(\Qq_S).$$

For $S=\{\infty\}$ and $D=1$, the space $X(1)^\infty$ is a covering of the level 1 modular curve 
$$X=\SL_2(\Zz)\bash \Hh\simeq \rmG(\Zz)\bash \rmG(\Rr)/K_\infty\simeq X(1)^\infty/K_\infty$$
with compact fibers.  More generally the $S$-arithmetic quotient $X(D)^{S}$ is a covering of the usual modular curve  $X_0(D)$ of level $D$. 
The compact fiber of a point $z \in \Hh$  is $$\Gamma(D)_Sg_zK_\infty\rmG(\Zz_S)$$ where $g_z\in\rmG(\Rr)$ is  such that $g_z.i=z$. 
We denote by
$$\pi_S : X(D)^S\ra X_0(D)\ra X_0(1)=X$$
the composite of this projection and the usual finite covering. Since the fibers of this map are compact, $\Gamma(D)_S\subseteq \rmG(\Qq_S)$ has finite covolume (i.e.\ is a lattice).

\subsection{Lifting the measures}To simplify notations, we omit the dependency in $D$ and write $X^S$ for $X(D)^S$, $\Gamma_S$ for $\Gamma(D)_S$ etc. 

Given a prime $q\not\in S$  and $a\in\Zz$ we set
$$u_{a/q} := \begin{pmatrix}
	1&a/q\\ 0&1
\end{pmatrix}.$$
We define
$$x^S_{q,a}:=\Gamma_S\Big(u_{a/q}\begin{pmatrix}
	1&0\\0&q
\end{pmatrix}	,u_{a/q},\cdots,u_{a/q}\Big)\in X^S.$$
We observe that since $q\not\in S$ the $\Gamma_S$-coset $x^{S}_{q,a}$ depends only on the congruence class $a\mods q$:  for $p\in S$ and $k\in \Zz$ we have 
$$u_{(a+qk)/q}=u_k.u_{a/q}\in \Gamma_Su_{a/q}.$$
Let 
$$H^{S}_{q}:=\{x^S_{q,a} \mid a\in\Zz/q\Zz\}\subseteq X^S.$$
By construction we have
$$\pi_S(H^{S}_{q})=H_q.$$
	Let us pass to the product: given $b\in (\Zz/q\Zz)^\times$ we set
$$H^{S}_{q,b}:=\big\{(x^S_{q,a},x^S_{q,ab}) \mid a\in\Zz/q\Zz\big\}\subseteq  X^S\times X^S$$
so that 
$$(\pi_S\times\pi_S)(H^{S}_{q,b})=H_{q,b}\in X\times X.$$
Let us denote by 
$\mu_{q,b}^{S}$ the uniform probability measure on $X^S\times X^S$ supported on $H^{S}_{q,b}$.

There are multiple advantages of lifting the whole situation to the $S$-arithmetic group quotient 
$X^S=\Gamma_S\bash \rmG(\Qq_S)$.  

\begin{itemize}
    \item The space $X^S$ is endowed with an action of $\rmG(\Qq_S)$ by right multiplication (noted $x.g=\Gamma_Sug$ for $x=\Gamma_Su\in X^S$ and $g\in \rmG(\Qq_S)$) and a left action on its space of functions: for $f:X^S\ra \Cc$ 
$$g.f:x\in X^S\mapsto f(x.g)\in\Cc.$$

\item We can test the lifted measures $\mu_{q}^{S}$ and $\mu_{q,b}^{S}$ against more general functions, e.g.\ corresponding to automorphic  forms not necessarily of weight $0$.  

\item The lifted measures have additional invariance from the places in $S$.
\end{itemize}

\subsection{Properties of the measures $\mu_{q,b}^{S}$}

We now specify the set of places $S$ we will use: given $D\geq 1$, let $q_1,q_2$ be two fixed primes coprime with $D$ and put 
$$S=\{\infty,q_1, q_2\}.$$
In this section we verify that any weak-$\star$ limit of the
$\mu_{q,b}^{S}$ is a joining. Let
$\mu^{S}_{q,b,i}, i=1,2$, denote the image of $\mu_{q,b}^{S}$ on $X^S$ under the first and second coordinates projections $\pi_i:X^S\times X^S\rightarrow X^S.$

\begin{lemma} As $q\ra\infty$, the measures $\mu^{S}_{q,b,i}$, $i=1,2$, converge to the Haar probability measure $\mu^S_{\rmG}$ on $X^S$.
\end{lemma}

\textbf{Proof.} For $q$ large enough, we have $q\not\in S$ and $(q,D)=1$, and the two projections of $H^S_{q,b}$ are given by the set $H_{q}^S$. Since $q$ is prime, $H_q^S$ is the $q$-Hecke orbit of the identity class $\Gamma_S\mathrm{Id}\in X^S$ minus one point, namely $x_{q,\infty}^{S}$ defined like $x^S_{q,a}$ but with with $u_{a/q}$ replaced by $\text{diag}(q,q^{-1})$;  equidistribution $H_{q}^S$ follows from the equidistribution of $q$-Hecke orbit.
\qed

\begin{remark}
	In particular any weak-$\star$ limit of the
$(\mu_{q,b}^{S})_{q}$ (the limit of a converging subsequence) is a probability measure.
\end{remark}
We now make use of the auxiliary places  $q_1$ and $q_2$. We introduce the two diagonal matrices
$$t_1=\begin{pmatrix}
	q_1^{-1}&\\&q_1
\end{pmatrix},\ t_2=\begin{pmatrix}
	q_2^{-1}&\\&q_2
\end{pmatrix}$$
which we view as embedded diagonally  into $\rmG(\Qq_S)$ 
and also  diagonally embedded into $\rmG(\Qq_S)\times \rmG(\Qq_S)$, i.e.\ via 
$t_1^\Delta=(t_1,t_1), t_2^\Delta=(t_2,t_2)\in \rmG(\Qq_S)\times \rmG(\Qq_S).$
\begin{lemma} The set $H^S_{q,b}$ is invariant under right multiplication by the  elements $t_1^\Delta,t_2^\Delta$.
\end{lemma}
\textbf{Proof.} For any $a\in\Zz/q\Zz$ we have
\begin{eqnarray*}
	x^S_{q,a}.t_1&=& \Gamma_S\Big(u_{a/q}\begin{pmatrix}
	1&0\\0&q
\end{pmatrix}\begin{pmatrix}
	q_1^{-1}&\\&q_1
\end{pmatrix}	,u_{a/q}\begin{pmatrix}
	q_1^{-1}&\\&q_1
\end{pmatrix},\cdots,u_{a/q}\begin{pmatrix}
	q_1^{-1}&\\&q_1
\end{pmatrix}\Big)\\
&=&\Gamma_S t_1\Big(u_{aq_1^2/q}\begin{pmatrix}
	1&0\\0&q
\end{pmatrix},u_{aq_1^2/q},\cdots,u_{aq_1^2/q}\Big)
\\
&=&\Gamma_S\Big(u_{aq_1^2/q}\begin{pmatrix}
	1&0\\0&q
\end{pmatrix},u_{aq_1^2/q},\cdots,u_{aq_1^2/q}\Big)=x^{S}_{q,aq_1^2}
\end{eqnarray*}
where we  used that $\text{diag}(q_1^{-1},q_1)$ commutes with $\text{diag}(1, q
)$ and  $t_1\in \Gamma_S$.  Also recall that $q_1$ is coprime to $q$ (since $q$ is a sufficiently large prime). 
Likewise
$$x^S_{q,ab}.t_1=x^S_{q,aq_1^2b}$$ and so
$$H^S_{q,b}.t_1^{\Delta}=\big\{(x^S_{q,aq_1^2},x^S_{q,aq_1^2b}) \mid  a\in\Zz/q\Zz\big\}=H^S_{q,b}.$$
The same computation applies to $t_2$.
\qed\\

In $\rmG(\Rr\times\Qq_{q_1}\times\Qq_{q_2})$ we use the factorisation 
$$t_1=\biggl(\begin{pmatrix}
	q_1^{-1}&\\&q_1
\end{pmatrix},\begin{pmatrix}
	q_1^{-1}&\\&q_1
\end{pmatrix},\mathrm{Id}_2\biggr)\biggl(\mathrm{Id}_2,\mathrm{Id}_2,\begin{pmatrix}
	q_1^{-1}&\\&q_1
\end{pmatrix}\biggr)=t'_1k_1,$$
$$t_2=\biggl(\begin{pmatrix}
	q_2^{-1}&\\&q_2
\end{pmatrix},\mathrm{Id}_2,\begin{pmatrix}
	q_2^{-1}&\\&q_2
\end{pmatrix}\biggr)\biggl(\mathrm{Id}_2,\begin{pmatrix}
	q_2^{-1}&\\&q_2
\end{pmatrix},\mathrm{Id}_2\biggr)=t'_2k_2,$$
say, so that $t'_1,t'_2$ satisfy assumption $\mcA'$ of \cite{EL}*{Definition 1.3}, while $k_1,k_2$ are contained in a  compact subgroup of the group of diagonal matrices of determinant $1$. Applying now \cite{EL}*{Cor.\ 3.4} (see also \cite{Kh}*{Thm.\ 4.4} for the formulation we use here), we obtain
\begin{prop}\label{alternative} Let $\mu_\infty$ be a weak-$\star$ limit of the measures  $\mu^S_{q,b}$ parametrized by pairs $(b, q)$ with $q \rightarrow \infty$. Then $\mu_\infty$ is a $(t_1^\Delta,t_2^\Delta)$-invariant probability measure, and any of its ergodic components is of the shape
\begin{itemize}
    \item[--] either $\mu_{\rmG\times\rmG}$, the (image of the) product Haar measure $\mu_{\rmG(\Qq_S)\times \rmG(\Qq_S)}$ on $X^S\times X^S$;
    \item[--] or $\mu_{\rmG,h}$, the (image of the) Haar measure $\mu_{\rmG(\Qq_S)}$ on a $\rmG(\Qq_S)$-orbit of the shape
    \begin{equation*}
        (\Gamma_S\times\Gamma_S)\rmG^\Delta(\Qq_S)(1,h)\subset X^S\times X^S,\ h\in\rmG(\Qq_S),
    \end{equation*}
    where $\rmG^\Delta$ denotes  the image of the diagonal embedding
    $\Delta:g\in \rmG\hookrightarrow (g,g)\in \rmG\times\rmG.$ 
    \end{itemize}
    More precisely, there exists a probability measure $\lambda$ on $\rmG(\Qq_S)$ such that $\mu_\infty$ is a convex  combination of $\mu_{\rmG\times\rmG}$ and of  
    $$\mu_{\rmG^\Delta}:=\int_{\rmG(\Qq_S)}\mu_{\rmG,h}d\lambda(h),$$
    i.e.\ there exists $c\in [0,1]$ such that
    $$\mu_\infty=(1-c)\mu_{\rmG\times\rmG}+c\mu_{\rmG^\Delta}.$$
\end{prop}

\subsection{Excluding the diagonal measures}
In this section we use the bounds for certain well chosen Weyl sums established in Theorem \ref{thm0} to show that the coefficient $c$ of $\mu_{\rmG^\Delta}$ in the decomposition of is zero, effectively concluding the proof of Theorem \ref{main1}. The proof of Corollary \ref{maincor} is a consequence of this result and \cite{EL}*{Cor 1.5}.\\

Consider a subsequence of the $\mu^S_{q,b}$ converging to some measure of the shape $$\mu_\infty=(1-c)\mu_{\rmG\times\rmG}+c\mu_{\rmG^\Delta}$$ and suppose that $c>0$.

Let $\Pi$ be a cuspidal automorphic representation such that its space $\Pi^{K^{(S)}}$ of $K^{(S)}$-invariant vectors occurs in $X^S$, and let $\vphi_1, \vphi_2$ be two smooth $L^2$-normalized vectors in $\Pi$.  
Consider the function on $X^S\times X^S$ given by
$$\vphi=\vphi_1\otimes\ov {\vphi}_2:(g,g')\in X^S\times X^S\mapsto \vphi_1(g)\ov{\vphi}_2(g').$$
Then for $h\in G(\Qq_S)$ we have 
$$\mu_{G,h}(\vphi)=\int_{X^S}\vphi_1(g)\ov\vphi_2(g.h)dg=\peter{\vphi_1,h.\vphi_2}$$
where we have set 
$$\peter{\vphi_1,\vphi_2}=\mu_G^S(\vphi_1\ov \vphi_2)=\int_{X^S}\vphi_1(g)\ov \vphi_2(g)dg.$$
We have therefore
$$\mu_{\rmG^\Delta}(\vphi)=\int_{\rmG(\Qq_S)}\peter{\vphi_1,h.\vphi_2}d\lambda(h)=\peter{\vphi_1,\lambda*\vphi_2}$$
where
$$\lambda* \vphi_2:g\mapsto \int_{\rmG(\Qq_S)} \vphi_2(g.h)d\lambda(h)$$
The integral $\lambda* \vphi_2$ converges  since the measure $
\lambda$ is finite and cusp forms have rapid decay. Roughly speaking, we will choose $\vphi_1 = \lambda \ast \vphi_2$. On the one hand, Theorem \ref{thm0} will imply that under certain conditions $\mu_{\infty}(\vphi) = 0$. On the other hand, we will show under certain conditions $\lambda \ast \vphi_2 \not= 0$, which is only possible if $c = 0$. We now make this precise. 

We make the following choices: let $K = \mathbb{Q}(\sqrt{229})$ (notice that 229 is a prime fundamental discriminant); this quadratic field has class number 3. Let $\chi$ be a class group character of order 3 ; it lifts to an automorphic form $f_{\chi}$ of level 229 with central character the Legendre symbol $\chi_{229} = (\frac{229}{.})$. Since $229 \equiv 1$ (mod 4), there exists  primitive Dirichlet character $\psi$ modulo 229 of order $4$ and its square is $\chi_{229}$. The automorphic form 
\begin{equation}\label{f}
   f = f_{\chi} \times \psi
   \end{equation} is an automorphic form of level $229^2$ with trivial central character and Laplace eigenvalue $1/4$ (i.e.\ spectral parameter $t = 0$). Let $$q_1 = 37, \quad q_2 = 53.$$ These are split primes in $K$ lying below principal ideals. Moreover, $q_1$ and $q_2$ are 4-th power residues modulo 229, so $\psi$ is trivial on $q_1, q_2$. We conclude that the automorphic representation $\Pi$ generated by $f$ satisfies the following properties:
\begin{itemize}
\item it is cuspidal;
\item it has trivial central character;
\item it satisfies the Ramanujan conjecture at all places;
\item it is ramified only at $229$ which is disjoint from $S$;
\item it has trivial Langlands parameters at all three places in $S = \{\infty, q_1, q_2\}$.
\end{itemize}
Of course these properties could have been obtained in a rather general way, the explicit construction above is only for illustrative purposes. 
The key point is the last property, for the following reason: We consider the representation $\Pi_S = \Pi_{\infty} \otimes \Pi_{q_1} \otimes \Pi_{q_2}$. Since $\Pi$ is locally unramified principal series at all places in $S$, we can consider its induced model $\mathcal{I}_S$ which is the space of functions $F : \rmG(\mathbb{Q}_S) \rightarrow \mathbb{C}$ that transform like
\begin{equation}\label{rule}
F(z n[x] a[y] g) = |y|_S^{1/2}   F(g), \quad z \in Z(\rmG), \quad n[x] = \left(\begin{matrix} 1 & x\\ & 1\end{matrix}\right), \quad a[y] = \text{diag}(y, 1)
\end{equation}
where we use that the Langlands parameters at places in $S$ are trivial. (Note that $y$ is only well-defined up to units in the non-archimedean case, but this is irrelevant for the transformation rule.) We have a $\rmG(\mathbb{Q}_S)$-equivariant vector space homomorphism $\iota: \Pi_S \rightarrow \mathcal{I}_S$. 

The image of  the weight 0 Maa{\ss} form  $f$ defined in \eqref{f}, lifted to a vector in the adelic representation $\Pi$, projected onto $\Pi_S$ in the model $\mathcal{I}_S$,  corresponds to the function $F$ on $\rmG(\mathbb{Q}_S)$  that is 1 on ${\PSO}(2) \times \rmG(\mathbb{Z}_{q_1}) \times \rmG(\mathbb{Z}_{q_2})$ which is obviously non-negative on all of $\rmG(\mathbb{Q}_S)$ by \eqref{rule}. In particular, $\lambda \ast F$ is not the zero function, so there exists some $\delta > 0$ with 
$$0 < \delta < \| \lambda \ast \iota^{-1}F \|^2$$
where $\iota^{-1} F$ is an $L^2$-normalized spherical vector in $\Pi_S$ (or in $\Pi$). 

For $R > 0$ let $B_R \subseteq \GQS$ denote a ball of radius $R$ about the origin. Choose $R> 0$ so that 
$$\int_{B_R} d\lambda(g) \geq  1 - \delta,$$
and let us denote by $\lambda_R$   the measure restricted to $B_R$. We now choose $$\vphi_2 = \iota^{-1} F, \quad \vphi_1 = \lambda_R \ast \vphi_2.$$ Note that $\vphi_1$ is a smooth vector since $\lambda_R$ is compactly supported. Setting as above $\vphi = \vphi_1 \otimes \bar{\vphi}_2$, we have
$$\mu_{\infty}(\vphi) = (1-c)\mu_{G \times G}(\vphi) + c\mu_{G^{\Delta}}(\vphi) = c \mu_{G^{\Delta}}(\vphi) = c \langle  \lambda_R \ast \vphi_2, \lambda \ast \vphi_2 \rangle 
    $$
    Since $\| \vphi_2 \|^2 = 1$ and $\lambda$ is a probability measure, we obtain by Cauchy-Schwarz and the definition of $\lambda_R$ that 
    $$\big|\langle  \lambda_R \ast \vphi_2, \lambda \ast \vphi_2 \rangle - \langle  \lambda \ast \vphi_2, \lambda \ast \vphi_2 \rangle\big| \leq \delta, $$
    so that 
\begin{equation}\label{delta}
\mu_{\infty}(\vphi) \geq c( \| \lambda \ast \vphi_2 \|^2 - \delta).
\end{equation}

On the other hand, setting $g_{q,a}=u_{a/q}\text{diag}(1, q)$  
we have 
 \begin{equation}\label{zeroconverge}
    \mu^S_{q,b}(\vphi_1\otimes\ov \vphi_2)=\frac{1}{q}\sum_{a\mods q} \vphi_1(g_{q,a} )\ov \vphi_2(g_{q,ba} ) \rightarrow \mu_\infty(\vphi)
\end{equation}
as $q \rightarrow \infty$ (over primes).  The left hand side equals
$$\int_{\| g \| \leq R} \frac{1}{q} \sum_{a\, (\text{mod } q)} \vphi_2(g_{q, a} g) \overline{\vphi_2(g_{q, ba})}\, d\lambda(g).$$
By strong approximation, we can translate this back into classical language. Clearly, $$\vphi_2(g_{q, ba}) = f\Big(\frac{ba + i}{q}\Big).$$
Moreover, any $g\in G(\Qq_S)$ with $\|g\|\leq R$ can be written
$g = g_{\infty} \times g_{q_1} \times g_{q_2}$ with $$g_{\infty} \in \left(\begin{matrix} y_0 & x_0\\ 0 & 1\end{matrix}\right)K_\infty, \quad g_{q_j} \in \left(\begin{matrix} q_j^{r_j} & \xi_j\\ 0 & 1\end{matrix}\right){\rm G}(\mathbb{Z}_{q_j})$$ for some $y_0\in \mathbb{R}^{\times}$, $x_0 \in \mathbb{R}$, $r_j \in \mathbb{Z}$, $\xi_j \in \mathbb{Z}[q_j^{-1}]$,  then  
$$\vphi_2(g_{q, a} g) = f\left( \left(\begin{matrix} Y & r_0\\ 0 & 1\end{matrix}\right)\left(\begin{matrix} 1/q & a/q\\ & 1\end{matrix}\right)\left(\begin{matrix} y_0 & x_0\\ & 1\end{matrix}\right).i\right) $$
for some $r_0 \in \mathbb{Z}[q_1^{-1}, q_2^{-1}]$ and $Y = q_1^{-r_1}q_2^{-r_2}$. Matrix multiplication shows that this equals
$$f\Big(\frac{aY+  Yx_0 +i y_0 Y}{q} + r_0 \Big).$$
Thus we see that with the notation \eqref{defW}
$$\frac{1}{q} \sum_{a\, (\text{mod } q)} \vphi_2(g_{q, a} g) \overline{\vphi_2(g_{q, ba})} = \mathcal{W}_{\bar{f}, f} (\bar{b}Y, q; x_0Y, y_0Y, r_0).$$

The entries $x_0Y, y_0Y, r_0, Y$ are bounded in terms of $R$ and as $q \rightarrow \infty$ (through primes), we have $$s(q;\bar{b}Y) \asymp_Y s(q;b) \rightarrow \infty,$$
by our assumption. By Theorem \ref{thm0} and Lemma \ref{min} we obtain that 
$$\mathcal{W}_{\bar{f}, f} (\bar{b}Y, q; x_0Y, y_0Y, r_0)\ra 0,$$
uniformly for $\|g\|\leq R$ and therefore \eqref{zeroconverge} tends to zero. 
 
 This contradicts \eqref{delta} unless $c = 0$ and this completes the proof of Theorem \ref{main1}.\\

The proof of Theorem \ref{main2} is similar, based on Theorem \ref{thm-cont} instead which we need to apply with $yr_0$ instead of $y$. If $y =a/q + O(1/qQ)$, then $yr_0 = ar_0/q + O(r_0/qQ)$, so as long as $r_0$ is fixed, this changes (as in the previous proof) only the implicit constants.\\

\textbf{Acknowledgement:} We would like to thank Edgar Assing, Manfred Einsiedler, Elon Lindenstrauss, Manuel L\"uthi, Asbjørn Nordentoft, Dinakar Ramakrishnan, Peter Sarnak and Radu Toma, for very useful conversations and comments. We would like to thank the referees for their very careful reading of the manuscript and for pointing out some inaccuracies in a previous version.


\begin{thebibliography}{ALMW22}

\bib{AkaBBKI}{article}{
   author={Aka, Menny},
   title={Joinings classification and applications [after Einsiedler and Lindenstrauss]},
   Journal={Seminaire Bourbaki},
   volume={1185},
    year={2021},   
}

\bib{AES}{article}{
   author={Aka, Menny},
   author={Einsiedler, Manfred},
   author={Shapira, Uri},
   title={Integer points on spheres and their orthogonal lattices},
   journal={Invent. Math.},
   volume={206},
   date={2016},
   pages={379--396},
   issn={0020-9910},
}

\bib{AEW}{article}{
   author={Aka, Menny},
   author={Einsiedler, Manfred},
   author={Wieser, Andreas},
   title={Planes in four-space and four associated CM points},
   journal={Duke Math. J.},
   volume={171},
   date={2022},
   pages={1469--1529},
   issn={0012-7094},
}
		
\bib{ALMW}{article}{
   author={Aka, Menny},
   author={Luethi, Manuel},
   author={Michel, Philippe},
   author={Wieser, Andreas},
   title={Simultaneous supersingular reductions of CM elliptic curves},
   journal={J. Reine Angew. Math.},
   volume={786},
   date={2022},
   pages={1--43},
   issn={0075-4102},
}

\bib{Ba}{article}{
   author={Banks, William D.},
   title={Twisted symmetric-square $L$-functions and the nonexistence of
   Siegel zeros on ${\rm GL}(3)$},
   journal={Duke Math. J.},
   volume={87},
   date={1997},
   number={2},
   pages={343--353},
   issn={0012-7094},
   doi={10.1215/S0012-7094-97-08713-5},
}

\bib{Bl}{article}{
   author={Blomer, Valentin},
   title={Shifted convolution sums and subconvexity bounds for automorphic
   $L$-functions},
   journal={Int. Math. Res. Not.},
   date={2004},
   number={73},
   pages={3905--3926},
   issn={1073-7928},
  
}



\bib{BB}{article}{
   author={Blomer, Valentin},
   author={Brumley, Farrell},
   title={Simultaneous equidistribution of toric periods and fractional moments of L-functions},
   journal={to appear in J. Eur. Math. Soc.},
 
}

\bib{BBK}{article}{
 author={Blomer, Valentin},
   author={Brumley, Farrell},
   author={Khayutin, Ilya},
   title={The mixing conjecture under GRH},
   note={https://arxiv.org/abs/2212.06280},
 
}

\bib{BFKMM}{article}{
   author={Blomer, Valentin},
   author={Fouvry, \'{E}tienne},
   author={Kowalski, Emmanuel},
   author={Michel, Philippe},
   author={Mili\'{c}evi\'{c}, Djordje},
   title={On moments of twisted $L$-functions},
   journal={Amer. J. Math.},
   volume={139},
   date={2017},
   pages={707--768},
   issn={0002-9327},
}

\bib{BH}{article}{
   author={Blomer, Valentin},
   author={Harcos, Gergely},
   title={The spectral decomposition of shifted convolution sums},
   journal={Duke Math. J.},
   volume={144},
   date={2008},
   pages={321--339},
   issn={0012-7094},
}

 \bib{BM}{article}{
   author={Bruggeman, Roelof W.},
   author={Motohashi, Yoichi},
   title={A new approach to the spectral theory of the fourth moment of the
   Riemann zeta-function},
   journal={J. Reine Angew. Math.},
   volume={579},
   date={2005},
   pages={75--114},
   issn={0075-4102},
}

 

\bib{BG}{article}{
   author={Bump, Daniel},
   author={Ginzburg, David},
   title={Symmetric square $L$-functions on ${\rm GL}(r)$},
   journal={Ann. of Math. (2)},
   volume={136},
   date={1992},
   number={1},
   pages={137--205},
   issn={0003-486X},
   doi={10.2307/2946548},
}

  \bib{CHT}{article}{
   author={Clozel, Laurent},
   author={Harris, Michael},
   author={Taylor, Richard},
   title={Automorphy for some $l$-adic lifts of automorphic mod $l$ Galois
   representations},
   journal={Publ. Math. Inst. Hautes \'{E}tudes Sci.},
   number={108},
   date={2008},
   pages={1--181},
   issn={0073-8301},
}

 

\bib{Du}{article}{
   author={Duke, W.},
   title={Hyperbolic distribution problems and half-integral weight Maass
   forms},
   journal={Invent. Math.},
   volume={92},
   date={1988},
   pages={73--90},
   issn={0020-9910},
}


\bib{DFI}{article}{
    AUTHOR = {Duke, W.}, 
    author={Friedlander, J.},
    author={Iwaniec, H.},
     TITLE = {Bounds for automorphic {$L$}-functions},
   JOURNAL = {Invent. Math.},
    VOLUME = {112},
      YEAR = {1993},
     PAGES = {1--8},
      ISSN = {0020-9910},
   MRCLASS = {11F66 (11F11)},
       URL = {https://doi.org/10.1007/BF01232422},
}
		
 

\bib{EL}{article}{
   author={Einsiedler, Manfred},
   author={Lindenstrauss, Elon},
   title={Joinings of higher rank torus actions on homogeneous spaces},
   journal={Publ. Math. Inst. Hautes \'{E}tudes Sci.},
   volume={129},
   date={2019},
   pages={83--127},
   issn={0073-8301},
}

\bib{ELMV}{article}{
   author={Einsiedler, Manfred},
   author={Lindenstrauss, Elon},
   author={Michel, Philippe},
   author={Venkatesh, Akshay},
   title={Distribution of periodic torus orbits and Duke's theorem for cubic
   fields},
   journal={Ann. of Math. (2)},
   volume={173},
   date={2011},
   pages={815--885},
   issn={0003-486X},
}
		



\bib{EMV}{article}{
   author={Ellenberg, Jordan S.},
   author={Michel, Philippe},
   author={Venkatesh, Akshay},
   title={Linnik's ergodic method and the distribution of integer points on
   spheres},
   conference={
      title={Automorphic representations and $L$-functions},
   },
   book={
      series={Tata Inst. Fundam. Res. Stud. Math.},
      volume={22},
   },
   date={2013},
   pages={119--185},
}

\bib{ELM}{article}{
   author={Einsiedler, Manfred},
   author={Lindenstrauss, Elon},
   author={Michel, Philippe},
   
   title={Arithmetic joint equidistribution},
   journal={unpublished manuscript},
   
   date={2018},
}
		



\bib{EMS}{article}{
   author={Elliott, P. D. T. A.},
   author={Moreno, C. J.},
   author={Shahidi, F.},
   title={On the absolute value of Ramanujan's $\tau $-function},
   journal={Math. Ann.},
   volume={266},
   date={1984},
   pages={507--511},
   issn={0025-5831},
}





\bib{Er}{article}{
   author={Erd\"{o}s, P.},
   title={On the sum $\sum^x_{k=1} d(f(k))$},
   journal={J. London Math. Soc.},
   volume={27},
   date={1952},
   pages={7--15},
   issn={0024-6107},
}

\bib{HST}{article}{
   author={Harris, Michael},
   author={Shepherd-Barron, Nick},
   author={Taylor, Richard},
   title={A family of Calabi-Yau varieties and potential automorphy},
   journal={Ann. of Math. (2)},
   volume={171},
   date={2010},
   number={2},
   pages={779--813},
   issn={0003-486X},
}

\bib{HL}{article}{
   author={Hoffstein, Jeffrey},
   author={Lockhart, Paul},
   title={Coefficients of Maass forms and the Siegel zero},
   note={With an appendix by Dorian Goldfeld, Hoffstein and Daniel Lieman},
   journal={Ann. of Math. (2)},
   volume={140},
   date={1994},
   number={1},
   pages={161--181},
   issn={0003-486X},
   doi={10.2307/2118543},
}

\bib{Ho}{article}{
   author={Holowinsky, Roman},
   title={Sieving for mass equidistribution},
   journal={Ann. of Math. (2)},
   volume={172},
   date={2010},
   pages={1499--1516},
   issn={0003-486X},
}



\bib{Kh}{article}{
   author={Khayutin, Ilya},
   title={Joint equidistribution of CM points},
   journal={Ann. of Math. (2)},
   volume={189},
   date={2019},
   pages={145--276},
   issn={0003-486X},
}

\bib{Kimsh}{article}{
   author={Kim, Henry H.},
   author={Shahidi, Freydoon},
   title={Cuspidality of symmetric powers with applications},
   journal={Duke Math. J.},
   volume={112},
   date={2002},
   number={1},
   pages={177--197},
   issn={0012-7094},
   doi={10.1215/S0012-9074-02-11215-0},
}


\bib{KimSar}{article}{
   author={Kim, Henry H.},
   title={Functoriality for the exterior square of ${\rm GL}_4$ and the
   symmetric fourth of ${\rm GL}_2$},
   note={With appendix 1 by Dinakar Ramakrishnan and appendix 2 by Kim and
   Peter Sarnak},
   journal={J. Amer. Math. Soc.},
   volume={16},
   date={2003},
   pages={139--183},
   issn={0894-0347},
}
		

\bib{LMW}{article}{
   author={Lindenstrauss, Elon},
   author={Mohammadi,A. },
   author={Wang, Z,},
   title={Effective equidistribution results for some one-parameter unipotent flows},
   note={https://arxiv.org/abs/2211.11099},
   }

\bib{MV}{article}{
   author={Michel, Philippe},
   author={Venkatesh, Akshay},
   title={Equidistribution, $L$-functions and ergodic theory: on some
   problems of Yu. Linnik},
   conference={
      title={International Congress of Mathematicians. Vol. II},
   },
   book={
      publisher={Eur. Math. Soc., Z\"{u}rich},
   },
   date={2006},
   pages={421--457},
}

\bib{NT}{article}{
   author={Nair, Mohan},
   author={Tenenbaum, G\'{e}rald},
   title={Short sums of certain arithmetic functions},
   journal={Acta Math.},
   volume={180},
   date={1998},
   pages={119--144},
   issn={0001-5962},
}

\bib{NewTh}{article}{
   author={Newton, James},
   author={Thorne, Jack A.},
   title={Symmetric power functoriality for holomorphic modular forms},
   journal={Publ. Math. Inst. Hautes \'{E}tudes Sci.},
   volume={134},
   date={2021},
   pages={1--116, 117--152},
    title={Symmetric power functoriality for holomorphic modular forms, I \& II},
   issn={0073-8301},
}

 



\bib{RW}{article}{
   author={Ramakrishnan, Dinakar},
   author={Wang, Song},
   title={On the exceptional zeros of Rankin-Selberg $L$-functions},
   journal={Compositio Math.},
   volume={135},
   date={2003},
   pages={211--244},
   issn={0010-437X},
}

\bib{Sa}{article}{
   author={Sarnak, Peter},
   title={Asymptotic behavior of periodic orbits of the horocycle flow and
   Eisenstein series},
   journal={Comm. Pure Appl. Math.},
   volume={34},
   date={1981},
   pages={719--739},
   issn={0010-3640},
}

\bib{SU}{article}{
    AUTHOR = {Sarnak, Peter},
    author =  {Ubis, Adri\'{a}n},
     TITLE = {The horocycle flow at prime times},
   JOURNAL = {J. Math. Pures Appl. (9)},
    VOLUME = {103},
      YEAR = {2015},
     PAGES = {575--618},
      ISSN = {0021-7824},
}
		



\bib{Saw}{article}{
   author={Sawin, Will},
   title={Bounds for the stalks of perverse sheaves in characteristic $p$
   and a conjecture of Shende and Tsimerman},
   journal={Invent. Math.},
   volume={224},
   date={2021},
   pages={1--32},
   issn={0020-9910},
}

\bib{Ta}{article}{
   author={Taylor, Richard},
   title={Automorphy for some $l$-adic lifts of automorphic mod $l$ Galois
   representations. II},
   number={108},
   date={2008},
   pages={183--239},
   issn={0073-8301},
}

\bib{Strom}{article}{
   author={Str\"{o}mbergsson, Andreas},
   title={On the uniform equidistribution of long closed horocycles},
   journal={Duke Math. J.},
   volume={123},
   date={2004},
   pages={507--547},
   issn={0012-7094},
}

\bib{Tak}{article}{
   author={Takeda, Shuichiro},
   title={The twisted symmetric square $L$-function of $\mathrm{GL}(r)$},
   journal={Duke Math. J.},
   volume={163},
   date={2014},
   number={1},
   pages={175--266},
   issn={0012-7094},
   doi={10.1215/00127094-2405497}
}

\bib{Tsi}{article}{
   author={Shende, Vivek},
   author={Tsimerman, Jacob},
   title={Equidistribution in ${\rm Bun}_2(\mathbb{P}^1)$},
   journal={Duke Math. J.},
   volume={166},
   date={2017},
   pages={3461--3504},
   issn={0012-7094},
}

\bib{Ven}{article}{
   author={Venkatesh, Akshay},
   title={Sparse equidistribution problems, period bounds and subconvexity},
   journal={Ann. of Math. (2)},
   volume={172},
   date={2010},
   pages={989--1094},
   issn={0003-486X},
}

\bib{Wo}{article}{
   author={Wolke, Dieter},
   title={Multiplikative Funktionen auf schnell wachsenden Folgen},
   journal={J. Reine Angew. Math.},
   volume={251},
   date={1971},
   pages={54--67},
   issn={0075-4102},
}

    
\end{thebibliography}
\end{document}